\documentclass[11pt]{amsart}
\numberwithin{equation}{section}
\newtheorem{theorem}{Theorem}
\newtheorem{definition}[theorem]{Definition}
\newtheorem{proposition}[theorem]{Proposition}
\newtheorem{lemma}[theorem]{Lemma}
\newtheorem{corollary}[theorem]{Corollary}
\numberwithin{theorem}{section}
\newtheorem{remark}[theorem]{Remark}
\def\R{\bf R}

\def\al{\aligned}
\def\eal{\endaligned}
\def\M{{\bf M}}
\def\be{\begin{equation}}
\def\ee{\end{equation}}
\def\lab{\label}
\def\a{\alpha}
\def\b{\beta}
\def\e{\epsilon}
\def\t{\tilde}
\def\R{\bf R}
\def\lam{\lambda}

\def\al{\aligned}
\numberwithin{equation}{section}

\begin{document}

\tracingpages 1
\title[no breathers theorem]{\bf A No breathers theorem for some noncompact Ricci flows}
\author{Qi S. Zhang}
\address{Department of
Mathematics Nanjing University, Nanjing 210093, China;
Department of
Mathematics, University of California, Riverside, CA 92521, USA }
\date{May  2012, revised May 2013}

\begin{abstract}
Under suitable conditions near infinity and assuming boundedness
of curvature tensor, we prove a no breathers theorem in the spirit
of Ivey-Perelman for some noncompact Ricci flows. These include
Ricci flows on asymptotically flat (AF) manifolds with positive
scalar curvature, which was studied in \cite{DM:1} and \cite{OW:1}
in connection with general relativity.  Since the method for the
compact case faces a difficulty, the proof involves solving a new
non-local elliptic equation which is the Euler-Lagrange equation
of a scaling invariant log Sobolev inequality.

It is also shown that the
 Ricci flow on AF manifolds with positive scalar curvature is uniformly $\kappa$ noncollapsed
 for all time.  This result, being different from Perelman's local noncollapsing result which holds in finite time,
 seems to have implications for the issue of longtime convergence.

\end{abstract}
\maketitle
\section{Statement of result}

A basic question in the study of the Ricci flow is: Are periodic orbits called  breathers trivial?  Here triviality means that metrics only move by diffeomorphisms and scaling
through out the period.
A Ricci flow  $(M, g(t))$, $t \in [t_1, t_2]$,  is called a breather if
 there is a positive constant
 $c$ and a diffeomorphism $\Psi$ on $M$ so that $ g(t_2)= c \Psi^* (g(t_1))$.
Perelman's no breathers
 theorems (\cite{P:1} Sections 2, 3) say that all periodical solutions of compact Ricci flows are gradient
 Ricci solitons, and hence trivial. See also earlier proof of this result by Ivey \cite{I:1} in three dimension steady and expanding case,
and \cite{Ca:1} and \cite{L:1} for further development on compact breathers. However, similar
result in the noncompact case is conspicuously absent. Finding nontrivial periodic orbits
has always been an useful topic in the  study evolution equations, which also include Ricci flows.
As indicated in the paper
\cite{OSW:1}, the nonexistence
of nontrivial breathers is associated to the irreversibility of world sheet in renormalization
group flow in string theory.  See also the papers \cite{FLW:1} and \cite{AKW:1} for further
motivations coming from physics, where the authors wish to rule of solitons which are also breathers by definition. Ruling out nontrivial breathers is also helpful in the
study of long time convergence problem in Ricci flow. For example, suppose one knows that
a Ricci flow $(M, g(t))$, $t \in [k, k+1]$, $k \to \infty$,  converges in $C^\infty_{loc}$ sense to a limit Ricci flow $(M_\infty, g_\infty(s))$, $s \in [0, 1]$. If the end points $(M_\infty, g_\infty(0))$ and $(M_\infty, g_\infty(1))$ differ only by scaling and diffeomorphism, then a no breather theorem
would imply that $(M_\infty, g_\infty(s))$, $s \in [0, 1]$ is a gradient Ricci soliton.
Actually  Theorem \ref{thnobreather} below implicitly implies that if certain scaling
invariant log Sobolev functionals at the end points share the same infimum which can
be reached by a minimizer, then $(M_\infty, g_\infty(s))$, $s \in [0, 1]$ is a gradient Ricci soliton.
This condition on the log Sobolev functional can be verified for many manifolds, including
the asymptotically flat ones. See Corollary \ref{coAF} below.

The purpose of this paper is to prove a no breathers theorem for some noncompact Ricci flows.
  Some
times an extension of a theorem from the compact case to a
noncompact one merely involves some technical improvements of the
method, plus some extra conditions near infinity. However the no
breathers theorem is different for two reasons. First, noncompact
Ricci flows arise naturally as the blow up limits of finite time
singularity of compact Ricci flows. In fact, most of the essential
singularity models for compact Ricci flows are noncompact. This
includes the well known cylinder $S^2 \times \R$ in the 3
dimensional case, which is also a trivial breather. Thus even if one is only interested in compact
Ricci flows, one still needs to study noncompact Ricci flows.
Second, the method of proof by Perelman for the no breathers
theorem does not seem to work for the noncompact case, especially
for the steady breather case. Recall that Perelman introduces the
$F$ functional which is defined as $F(v) = \int_{M} ( 4 |\nabla v
|^2 + R v^2) dg$ where $R$ is the scalar curvature of the manifold
and $v \in W^{1, 2}(M)$ and $\Vert v \Vert_{L^2(M, g)} =1$. He
proved that the infimum of $F$ is a nondecreasing function of time
along a Ricci flow $(M, g(t))$; moreover it is a constant if and
only if the Ricci flow is a steady gradient soliton. Using the
fact that the infimum is reached by a minimizer when $M$ is
compact, Perelman proved that there is no nontrivial steady
breathers for compact Ricci flows, i.e. a steady breather is
necessarily a steady gradient soliton. If one attempts to extend
this argument to noncompact Ricci flow, one faces an immediate
difficulty. Namely, the infimum of the $F$ functional is not
reached by a function on a typical noncompact manifold such as
$\R^n$ or $S^2 \times \R$. In fact, on $\R^n$, the $F$ functional
is nothing but the Dirichlet energy (multiplied by $4$) and it is
well known that there is no $L^2$ minimizer. For this reason, we
need to look for a different method.

In this paper, we consider the functional (\ref{lvg}). When the parameter $\a = 1$, it is
the limiting case of Perelman's $W$ entropy and which can be
regarded as a scaling invariant version of the Log Sobolev
inequality introduced by Weissler \cite{W:1}. The corresponding
Euler-Lagrange equation is a nonlocal, nonlinear elliptic
equation. Unlike the $F$ functional, the minimizer of (\ref{lvg})
exists on many typical noncompact manifolds. Using this we prove
a no breathers theorem on some noncompact Ricci flows.
The study of the functional (\ref{lvg}) and its minimizer equation
potentially has further applications.

Let's introduce notations and definitions to be used in the paper.
We use $M$ to denote a $n (\ge 3)$ dimensional Riemannian manifold and $g(t)$ to denote the
metric at time $t$; $d(x, y, t)$ is the geodesic distance under $g(t)$; Unless stated otherwise, we
assume the curvature tensor is bounded at each time $t$.
$B(x, r, g(t)) = \{ y \in {\M} \ | \ d(x, y, t) < r \}$  is the geodesic ball of radius $r$, under metric $g(t)$,
centered at $x$, and $|B(x, r, t)|_{g(t)}$ is the volume of
$B(x, r, t)$ under $g(t)$;  when no confusion arises we may also use $B(x, r)$ or $B(x, r, t)$
to denote $B(x, r, g(t)) $;  $dg(t)$ is the volume element; $x_0$ is a reference point on $M$. We also
reserve $R=R(x, t)$ as the scalar curvature under $g(t)$.  A generic positive constant is denoted
by $C$ or $c$ whose value may change from line to line. When we say that a sequence of
pointed manifolds converges in $C^\infty_{loc}$ sense, we mean they converge in the usual
Cheeger-Gromov sense. That is, subject to diffeomorphisms, the metrics converge
in $C^\infty_{loc}$ sense.
The definition of asymptotically flat manifolds can be found in the beginning of Section 2.

\begin{definition} (Log Sobolev functionals, infimum, infimum at
infinity)
 \lab{deflv} Let $(M, g)$ be a $n$ dimensional Riemannian
manifold with metric $g$ and $D \subset M$ be a domain.

(a). Given  functions $ v \in W^{1, 2}_0(D, g)$ with $\Vert v
\Vert_{L^2(D)}=1$,  and a number $\alpha \ge 1$,  the log Sobolev
functionals with parameter $\a$ is defined by
\be
\lab{lvg}
\al
 L(v,
g, \a, D) &= - \int_D v^2 \ln v^2 dg + \a \frac{n}{2} \ln
\left(\int_D ( 4 |\nabla v |^2 + R v^2 ) dg + E^-_0 \right) + s_n \\
&\equiv - N(v) + \a  \frac{n}{2} \ln (F(v) + E^-_0) + s_n.
\eal
\ee
Here $R$ is the scalar curvature;  $E^-_0 = - \min \{0,  \inf \{ F(v) \, |, v \in C^\infty_0(D), \Vert v \Vert_{L^2}=1 \} \}$;
 $s_n=-\frac{n}{2} \ln (2 \pi n)
- \frac{n}{2}$.

(b). The infimum of the log Sobolev functional is denoted by
\[
\lam(g, \a, D) =
\inf \{ L(v, g, \a, D) \, | \,  v \in W^{1, 2}_0(D, g), \quad \Vert v \Vert_{L^2(D)}=1 \}.
\]

(c). When $\a = 1$ and $D=M$, the infimum of the log Sobolev
functional at infinity is
\[
\lam_\infty(g, 1, M) = \lim_{r \to \infty} \lam(g, 1, M-B(x_0, r))
\]where $x_0$ is a reference point in $M$.

\end{definition}

If $D = M$, then for simplicity we write
\[
L(v, g, \a)=L(v, g, \a, M), \qquad \lam(g, \a)=\lam(g, \a, M).
\]If $\a=1$, we may suppress $\a$ and write
\[
 L(v, g)=L(v, g, 1), \qquad \lam =\lam(g)=\lam(g, 1)=\lam(g, 1, M) \qquad
\lam_\infty = \lam_\infty(g)=\lam_\infty(g, 1, M).
\]

\begin{remark}  When $M= \bf R^n$ and $\a=1$, then $L(v, g)$ is the log
Sobolev functional introduced by Weissler \cite{W:1}, which is a
scaling invariant version of the log Sobolev functional originally
introduced by Gross \cite{G:1} and Federbush \cite{F:1}. Observe
that $\lam(g)$ is invariant under scaling and diffeomorphism. See
the beginning of proof of Theorem \ref{thnobreather} below.

$\lam(g)$ is related to Perelman's $\nu$ invariant in Section 3 of
\cite{P:1}. We are not sure if they are the same.

When $F(v)$ becomes $0$ but $L(v)$ is finite, the functional $L$ is
regarded as $-\infty$. When the scalar curvature $R \ge 0$, it is clear that $E^-_0=0$.
\end{remark}

\begin{definition} (gradient Ricci solitons)  A Riemannian manifold $(M, g)$ is called a
gradient Ricci soliton if there exists a smooth function $f$ on
$M$ and a constant $\e$ such that
\be
\lab{defGRSol}
Ric+ Hess f +
\frac{\e}{2} g =0.
 \ee

$(M, g)$ is called a expanding, steady and shrinking gradient Ricci soliton if
$\e>0, \e=0$ and $\e<0$ respectively.
\end{definition}

\medskip

The following is the main result of the paper.

\begin{theorem}
\lab{thnobreather}
 Let $(M, g(t))$,  $\partial_t g_{ij} = - 2 R_{ij}$, $t \in [0, T]$ be a complete, noncompact  Ricci flow
 with bounded curvature tensor and nonnegative scalar curvature.  Suppose $(M, g(t))$ is a breather, i.e.
 for two moments $t_1, t_2 \in [0, T]$, $t_1<t_2$, there is a positive constant
 $c$ such that $(M, c g(t_1))$ and $(M,  g(t_2))$ differ only by diffeomorphism.

 Suppose also the following conditions hold.

 (a) $-\infty<\lam(g(t_1)) < \lam_\infty(g(t_1))$.

 (b) Either $|B(x_0, r, t_1)|_{g(t_1)} \le C r^n$, for some $C>0$ and all $r>0$,
 or  $R(x, t_1) \ge \frac{C}{1+ d(x, x_0, t_1)^2}$ for some constant $C>0$.

 Then $(M, g(t))$ is a
  gradient Ricci soliton.
\end{theorem}
\medskip

\begin{remark}  Although Condition (a) looks similar to
a well known condition on the existence of point spectrum for the Laplacian on
 noncompact manifolds, however, our condition is much less restrictive in the case the scalar
curvature is nonnegative.  It is well known that the Laplacian on asymptotically flat (AF) manifolds  does not have a point spectrum. But  Proposition \ref{prAFsob}  (b) and
Proposition \ref{prlam<0} below say that  AF manifolds with positive scalar curvature
satisfy Condition (a) unless they are shrinking gradient solitons.

It would be interesting to find more manifolds such that Condition (a) holds. We suspect that
certain decay condition of the curvature near infinity is sufficient.
\end{remark}
\medskip

Naturally one is obliged to present some examples of Ricci flows
where the conditions of the theorem is met. Condition (a) is easy to be met  since one can modify the metric on a compact domain of a
manifold so that $\lam(g)$ becomes arbitrarily negative, while
$\lam_\infty(g)$ remains the same. Let $x_0$ be a reference point,
we can construct a metric $g(t_1)$ such that the volume of the
unit ball $B(x_0, 1)$ is very small but the scalar curvature is
bounded by $1$. A flat cylinder with small aperture is such an
example. So given a positive number $\kappa$, the manifold is
$\kappa$ collapsed at scale $1$. Hence $\lam(g(t_1))$ is very
negative. Indeed, by Proposition \ref{prlamtosob}, if
$\lam(g(t_1)) > - C > -\infty$, then $(M, g(t_1))$ is $\kappa$
non-collapsed below scale $1$. Here $C$ depends on $\kappa$. But
$\lam_\infty(g(t_1))$ is totally independent of $\lam(g(t_1))$.

Condition (b) is satisfied automatically by ancient $\kappa$
solutions of 3 dimensional Ricci flow, which include gradient
shrinking solitons with nonnegative sectional curvature. See \cite{P:1} and \cite{P:2}.

Another type of examples is the Ricci flow on asymptotically flat (AF) manifolds (c.f.
Definition \ref{defAF}), which is interesting due to connections
to general relativity. Useful properties of these kind of Ricci
flows have bee proven in \cite{DM:1}, \cite{OW:1}. For example,
they proved that the AF property is preserved under Ricci flow.

\begin{corollary}
\lab{coAF}
Let $(M, g(t))$ be a Ricci flow on an asymptotically flat  manifold with positive scalar curvature. If $(M, g(t))$ is a breather then it is a gradient Ricci soliton.
\end{corollary}
\proof By Proposition \ref{prAFsob}  (a), we know $\lam(M, g(t))>
-\infty$. If $(M, g(t))$ is a gradient shrinking Ricci soliton,
then the proof is done. So we assume $(M, g(t))$ is not a gradient
Ricci soliton.  By Theorem \ref{thallnoncol} and Proposition
\ref{prlamtosob}, $(M, g(t))$ is $\kappa$ noncollapsed.  Applying
Proposition \ref{prAFsob} (b) and Proposition \ref{prlam<0}, we
find that $\lam(M, g(t))<0 \le \lam_\infty(M, g(t))$. By
Definition of AF manifolds, we also have $|B(x_0, r, t)|_{g(t)}
\le C r^n$, for some $C=C(t)>0$. Therefore, all the conditions of
the theorem are satisfied and the conclusion follows. \qed

\medskip
{\remark  In a recent paper \cite{Ha:1}, Haslhofer considered
Ricci flows on some AF manifolds with positive scalar curvature.
Under the extra assumption that the scalar curvature is
integrable, he modified the domain of Perelman's $F$ entropy to
include only smooth functions converging to $1$ sufficiently fast
at infinity. Using the monotonicity of this modified $F$ entropy,
one can also prove that steady breathers are steady solitons in
this case, under further assumptions near infinity on the
diffeomorphism in the definition of breathers. Also a no breather
theorem for some noncompact Ricci flows in the case of shrinking
solitons is proven in \cite{Z:2}.}

\medskip

{\remark  One may wonder if a no breathers theorem still holds when the scalar curvature changes sign.
When the operator $-\Delta + R$ has a negative eigenvalue,  under mild assumptions near infinity,
one can prove that the eigenfunction decays to zero exponentially fast. Then one can use Perelman's
original method described earlier to prove that steady breathers are steady gradient solitons. However, steady gradient solitons are ancient solutions. According to
\cite{Ch:1}, the scalar curvature is nonnegative. So the operator $-\Delta + R$
can not have negative eigenvalue. This contradiction shows that  no steady breathers exist in this case.
}

Let us outline the proof of the theorem.  The main hurdle is to prove the following theorem
which states that
the infimum of the functional $L(v, g(t_2), 1, M)$ is reached by a smooth function in $W^{1, 2}(M, g(t_2))$.

\begin{theorem}
\lab{thminimizer} Let $(M, g)$ be a noncompact manifold with
bounded curvature and nonnegative scalar curvature, which also
satisfies

 (a) $-\infty<\lam(g) < \lam_\infty(g)$.

 (b) Either $|B(x_0, r)|_{g} \le C r^n$, for some $C>0$ and all $r>0$,
 or  $R(x) \ge \frac{C}{1+ d(x, x_0)^2}$ for some constant $C>0$.

Then there exists a minimizer $v$ for the Log Sobolev functional $L(\cdot, g, 1, M)$, which satisfies
the equation
\be
\lab{maineq}
 \frac{n}{2} \frac{4 \Delta v - R v}{ \int ( 4 | \nabla v |^2 + R v^2 ) dg }
+  2 v \ln v +   \left( \lam(g, 1,  M)+   \frac{n}{2} -  \frac{n}{2}  \ln  \int ( 4 | \nabla v |^2 + R v^2 ) dg
- s_n \right)  v = 0.
\ee
\end{theorem}

The proof is done by an approximation process that involves a
priori estimates and a blow up analysis. This strategy has been
used to study variational problems involving critical functionals.
Recently in \cite{DE:1}  Dolbeault and Esteban treated a similar
functional on the cylinder $S^n \times \R$. We benefitted from the
ideas in that paper. However, we are facing new difficulties since
our functional is scaling invariant and its component $\ln F(v)$
may not be bounded from below.  These make it difficult to apply
P. L. Lions' concentrated compactness method near infinity
directly. However  under the extra assumption $\lam(g(t_2))<
\lam_\infty(g(t_2))$, we can show that the Lions' method
\cite{Lio:1} works on special regions where the $L^2$ norm of $v$
has faster than usual decays. We also use a fact that a sequence
of Boltzmann  entropy $N(v_k)$ satisfies the reverse Fatou lemma
when $\{ v_k \}$ is a sequence of bounded functions with the same
$L^2$ norm. Once a minimizer is found, we can use Perelman's
monotonicity formula to show that $(M, g(t))$ is a gradient Ricci
soliton since $\lam(g(t_1))= \lam(g(t_2))$.

\section{preliminaries and all time $\kappa$ noncollapsing on AF manifolds}

In this section, we present  a number of elementary results to be
used in the proof of the theorems and the corollary. We also prove that  the
 Ricci flow on AF manifolds with positive scalar curvature is uniformly $\kappa$ noncollapsed
 for all time.

\begin{definition}
\lab{defAF}
A complete, noncompact Riemannian manifold $M$ is called Asymptotically Flat of order $\tau$ if
there is a partition $M=M_0 \cup M_\infty$, which satisfies the following properties.

(i). $M_0$ is compact.

(ii). $M_\infty$ is the disjoint union of finitely many components each of which is diffeomorphic to $({\R^n} - B(0, r_0))$ for some $r_0>0$.

(iii). Under the coordinates induced by the diffeomorphism,  the metric $g_{ij}$ satisfies, for $x \in
M_\infty$,
\[
g_{ij}(x) = \delta_{ij}(x) + O(|x|^{-\tau}), \quad \partial_k
g_{ij}(x) =  O(|x|^{-\tau-1}), \quad \partial_k  \partial_l
g_{ij}(x) =  O(|x|^{-\tau-2}).
\]
\end{definition}

\begin{remark}
\lab{re2.1}
For convenience we will equip the compact component $M_0$ with a reference point $0$.
We will also assume  that $M_\infty$ has only one connected
component. This assumption does not reduce any generality for Corollary \ref{coAF} and the results in this section.
Since the key inequality $\lambda_\infty(g) \ge  0$ always holds regardless the number of connected
components for $M_\infty$.
\end{remark}

According to Theorem (1.1) in \cite{BKN:1}, if $M$ has one end,  the curvature tensor decays sufficiently fast near infinity and
$|B(0, r)| \ge c r^n$ when $r$ is large, then $M$ is AF. Here $n$ is the dimension.

\begin{proposition}
\lab{prAFsob} Let $(M, g)$ be an AF manifold of dimension $n \ge
3$. Suppose the scalar curvature $R$ is positive everywhere.

(a). Then there exists a constant $A>0$,  such that
\be
\lab{alesob} \left( \int_{M} v^{2n/(n-2)} dg \right)^{(n-2)/n} \le
A \int_{M} ( 4 |\nabla v |^2 + R v^2 ) dg, \quad \forall  v \in
W^{1, 2}(M, g);
\ee
moreover $\lam(g)$ is bounded from below i.e.
\be \lab{alelogsob}
 \int_{M} v^2 \ln v^2 dg \le \frac{n}{2} \ln \left( A \int_{M} ( 4 |\nabla v |^2 + R v^2 ) dg
 \right),
 \ee $\forall v \in W^{1, 2}(M, g), \Vert v \Vert_{L^2(M, g)}=1.$

(b). $\lam_\infty(g) \ge 0$.
\end{proposition}

\proof

(a).  We just need to prove (\ref{alesob}) since (\ref{alelogsob})
follows from Jensen inequality.

Pick and fix $r_0>0$ sufficiently large, so that a coordinate
system on $M-B(x_0, r_0)$ exists, which satisfies the defining
condition of AF manifolds. Let $\phi \in C^\infty_0(M)$ be a
cut-off function such that $ \phi =1 $ on $B(0,  r_0)$, \quad $
\phi =0 $ on $M-B(0, 2 r_0)$, \quad $0 \le \phi \le 1$ and
$|\nabla \phi | \le C/r_0$. For any $v \in C^\infty_0(M)$, the
function $v (1-\phi)$ is supported in $M-B(0,  r_0)$.

Let $J : M-B(0, r_0) \to \R^n$ be the coordinate map. Then the function
\be
\lab{f=J^{-1}}
f \equiv [v (1-\phi)] \circ J^{-1}
\ee is a smooth, compactly supported function in $\R^n$, after extending by zero value.
By the Euclidean Sobolev inequality, the following inequality holds
\be
\lab{EucSob}
\left( \int_{\R^n} f^{2n/(n-2)} dx \right)^{(n-2)/n} \le S_0
\int_{\R^n} |\nabla_{\R^n} f |^2  dx \ee
where $dx$ is the
Euclidean volume element and $\nabla_{\R^n}$ is the Euclidean
gradient. According to the definition of AF manifolds, there
exists a positive constant $c$ such that
\be
\lab{dx=dg}
c^{-1} dx \le dg(x) \le c dx, \quad c^{-1} |\nabla_{\R^n} f | \le |\nabla [v (1-\phi)] | \le c
|\nabla_{\R^n} f |.
\ee Here $ |\nabla [v (1-\phi)] |$ is the length of the gradient of $v (1-\phi)$, both with respect to  $g$. Therefore, there exists a positive constant $C$ such that
\be
\lab{sobv(1-phi)}
\left( \int_M | v (1-\phi)^{2n/(n-2)} dg \right)^{(n-2)/n} \le
C \int_M  |\nabla [v (1-\phi)]^2  dg.
\ee By this and Minkowski inequality, together with the standard Sobolev inequality
in the ball $B(0, 2 r_0)$, we deduce
\be
\lab{sobv1}
\al
&\left( \int_M v^{2n/(n-2)} dg \right)^{(n-2)/n} \\
&\le  2 \left( \int_M | v (1-\phi) |^{2n/(n-2)} dg \right)^{(n-2)/n} + 2 \left( \int_M (v \phi) |^{2n/(n-2)} dg \right)^{(n-2)/n} \\
&\le
C \int_M  |\nabla [v (1-\phi)]^2  dg + C \int_M  |\nabla (v \phi)|^2  dg
+C \int_M  (v \phi)^2  dg.
\eal
\ee Hence
\be
\lab{sobv2}
\left( \int_M v^{2n/(n-2)} dg \right)^{(n-2)/n} \le C \int_M  |\nabla v  |^2 dg
+ C \sup |\nabla \phi |^2  \int_{B(0, 2 r_0)}  v^2  dg
+ C \int_{B(0, 2 r_0)}  (v \phi)^2  dg.
\ee Since $R(x)>0$ for every $x \in \R^n$ by assumption, this implies, for
some constant $A>0$, that
\be
\lab{sobv3}
\left( \int_M v^{2n/(n-2)} dg \right)^{(n-2)/n} \le A \int_M  (4 |\nabla v  |^2
+ R v^2) dg.
\ee This is (\ref{alesob}), i.e. part (a).

\medskip

Now we prove part (b).

\noindent First we prove the following assertion.

{\it When the radius $r$ is sufficiently large, we have
\be
\lab{assert2.1}
\lam(g, 1, M-B(0, r)) \ge \lam(g_E, 1, {\R}^n-J(B(0, r)) + o(1).
\ee Here $J$ is the coordinate map near infinity in the definition
of AF manifold; $o(1)$ is a quantity whose absolute value goes to
$0$ when $ r \to \infty$; $g_E$ is the Euclidean metric.}

Pick a function $v \in C^\infty_0(M-B(0, r))$ with $\Vert v
\Vert_{L^2}=1$. Given any $\e>0$, by definition of AF manifolds,
for $x \in M-B(0, r)$ with $r$ sufficiently large, there are the following relations
\be
\lab{2.100}
(1-\e) dx \le dg(x) =\sqrt{ det g(x)} dx \le (1+\e) dx,
\ee
\be
\lab{2.101}
(1-\e) |\nabla_{\R^n} f | \le |\nabla v | \le (1+\e)
|\nabla_{\R^n} f |
\ee where $f = v \circ J^{-1}$ and $J$ is the coordinate map. Also
$\nabla_{\R^n}$ is the Euclidean gradient.  Hence
\be
\lab{2.102}
\int_{M} ( 4 |\nabla v |^2 + R v^2 ) dg
\ge (1-\e )^2
 \int_{\R^n} 4 |\nabla_{\R^n} f|^2 \sqrt{ det g(x)} dx
\ee  Write $\sqrt{ det g(x)} = w^2$,  a routine calculation shows
\be
\lab{2.103}
\al
\int_{\R^n} 4 |\nabla_{\R^n} f|^2 \sqrt{ det g(x)} dx& = \int_{\R^n} 4 |\nabla_{\R^n} f|^2 w^2 dx\\
&=\int_{\R^n} 4 |\nabla_{\R^n}( f w) |^2 dx + 4 \int_{\R^n} (f w)^2 \frac{\Delta w}{w} dx.
\eal
\ee By definition of $AF$ manifolds, we know that $\frac{| \Delta w(x) |}{w(x)} \le \frac{c}{|x|^{2+\tau}}$ with
$\tau>0$. Hence, by the Hardy's inequality in the Euclidean space, we have
\be
\lab{2.104}
\int_{\R^n} 4 |\nabla_{\R^n} f|^2 \sqrt{ det g(x)} dx \ge (1+o(1)) \int_{\R^n} 4 |\nabla_{\R^n}( f w) |^2 dx,
\ee which implies
\be
\lab{fv>1-e}
\int_{M} ( 4 |\nabla v |^2 + R v^2 ) dg  \ge (1-\e)^2 (1+o(1)) \int_{\R^n} 4 |\nabla_{\R^n}( f w) |^2 dx.
\ee
Also
\be
\lab{2.105}
\al
\int_M v^2 \ln v^2 dg &=
 \int_{\R^n} (f w)^2 \ln f^2 dx =\int_{\R^n} (f w)^2 \ln (f w)^2 dx - \int_{\R^n} (f w)^2 \ln w^2 dx \\
 &=\int_{\R^n} (f w)^2 \ln (f w)^2 dx + o(1).
 \eal
\ee This and (\ref{fv>1-e}) imply
 that
\be
\lab{2.106}
L(v, g, 1, M-B(0, r)) \ge L(f w, g_E, 1, {\R}^n-J(B(0, r))) + o(1) - n \e.
\ee Since $\Vert f w \Vert_{L^2(\R^n)} = 1$, by taking the infimum of this inequality,
it is easy to see that
\be
\lab{2.107}
\lam(g, 1, M-B(0, r)) \ge \lam(g_E, 1,
{\R}^n-J(B(0, r)) + o(1) -n \e.
\ee Since $\e$ is arbitrary,  the assertion is proven.

Using $\lam(g_E, 1, {\R}^n-J(B(0, r)) \ge \lam(g_E, 1, {\R}^n) =
0$, we see that
\be
\lab{2.108}
\lam_\infty(g) = \lim_{r \to \infty} \lam(g, 1, M-B(0, r)) \ge 0.
\ee This proves part (b).
 \qed

\begin{proposition}
\lab{prlam<0}
Let $(M, g(t))$ be a noncompact Ricci flow on the time interval $(A, B)$ such that
the curvature tensor is bounded for each time $t \in (A, B)$. Suppose also $(M, g(t))$
is $\kappa$ noncollapsed below scale $1$ and the scalar curvature is nonnegative.
If $(M, g(t))$ is not a gradient shrinking soliton, then
\be
\lab{2.109}
\lam(g(t_0)) \equiv \lam(g(t_0), 1, M) < 0, \quad t_0 \in (A, B).
\ee Moreover, for any $x_0 \in M$, when $r_0$ is sufficiently large, we have
\be
\lab{2.110}
 \lam(g(t_0), 1,  B(x_0, r_0) )< 0.
\ee Here $B(x_0, r_0)  = B(x_0, r_0, g(t_0)) $.
\end{proposition}

\proof For compact Ricci flows, Perelman (\cite{P:1} Section 3) already proved a similar inequality for
his $\nu$ invariant. The following
proof for the noncompact case is similar, except that one needs to justify integration by parts near
infinity.

Without loss of generality we assume $t_0<0 \in (A, B)$.  Let $u=u(x, t)=G(x, t; x_0, 0)$ be the fundamental solution of the conjugate heat equation
\be
\lab{2.111}
\Delta u - R u + \partial_t u=0, \qquad t<t_0.
\ee Let $s = -t$ and
\be
\lab{went}
W(g(t), u(\cdot, t), t) =
\int_{M} \left[ s (\frac{|\nabla u |^2}{u} + R u ) - u \ln u -
\frac{n}{2} \ln (4 \pi s) u - n u \right] dg(t)
\ee be Perelman's $W$ entropy corresponding to $u=u(x, t)$.
According to \cite{P:1} Section 3,
\be
\lab{2.112}
\frac{d}{dt} W(g(t), u(\cdot, t), t) = 2 s \int | Ric_{g(t)} - Hess_{g(t)} \ln u - \frac{1}{2s}
 g(t) |^2 u dg(t)  \ge 0
\ee with strict inequality holding unless $(M, g(t))$ is a gradient shrinking soliton.
Moreover $\lim_{t \to 0} W(g(t), u(\cdot, t), t)=0$. We comment that Perelman proved the result
for compact Ricci flows. In the noncompact case one needs to justify the integrability of the
quantities involved. Since $(M, g(t))$ has bounded geometry within any finite time interval and
$u$, as fundamental solution has Gaussian decay near infinity, the integrability issue
has been worked out in \cite{cty2007} and \cite{chowetc3}  Chapter 19 e.g..

Since $(M, g(t))$ is not a gradient shrinking soliton, $
 \frac{d}{dt} W(g(t), u(\cdot, t), t)$ is strictly positive. From the assumption $t_0<0$, we obtain
\be
\lab{2.113}
W(g(t_0),  u(\cdot, t_0), t_0) < \lim_{t \to 0} W(g(t), u(\cdot, t), t) =0.
\ee Observe that with $\rho>0$ regarded as a free parameter and taking $v=\sqrt{u(\cdot, t_0)}$, we have
\be
\lab{2.114}
\al
L(\sqrt{u(\cdot, t_0)}, g(t_0), 1) &
=   - \int_M v^2 \ln v^2 dg(t_0) +  \frac{n}{2} \ln \left(\int_M ( 4 |\nabla v |^2 + R v^2 ) dg(t_0) \right)
+ s_n\\
&=\inf_{\rho >0} \int_{M} \left[ \rho (\frac{|\nabla u |^2}{u} + R u ) - u \ln u -
\frac{n}{2} \ln (4 \pi \rho) u - n u \right] dg(t_0) \\
&\le \int_{M} \left[ |t_0| (\frac{|\nabla u |^2}{u} + R u ) - u \ln u -
\frac{n}{2} \ln (4 \pi |t_0|) u - n u \right] dg(t_0)\\
&=W(g(t_0),  u(\cdot, t_0), t_0)<0.
\eal
\ee Here $u=u(\cdot, t_0)$ and $R=R(\cdot, x_0)$. This shows, since $\lam(g(t_0))$ is the
infimum of the log Sobolev functional $L$, that
$\lam(g(t_0), 1) < 0$.

The second statement of the lemma is an easy consequence of the fact that
$\lam(g(t_0)) = \lim_{r_0 \to \infty} \lam(g(t_0), 1, B(x_0, r_0)).$
 \qed

 \begin{proposition}
\lab{prlamal}
Let $(M, g)$ be a noncompact manifold  such that $\lam(g)> -\infty$.

(a). For any
$x_0 \in M$,  $r_0>0$, and for all $\a \ge 1$,
the infimum of the log Sobolev functionals $L(\cdot, g, \a, B(x_0, r_0))$ satisfy:
\[
\lam(g, \a, B(x_0, r_0) ) \ge -C
\]where $C$ is a constant depending only on $\a$, $n$, the constant $\lam(g)$  and $|B(x_0, r_0)|$.

(b). $\lim_{\a \to 1^{+}} \lam(g, \a, B(x_0, r_0) )= \lam(g, 1, B(x_0, r_0))$.
\end{proposition}

\proof For simplicity we  use $B$ to denote $B(x_0, r_0)$ and $E^-_0=0$ in the proof.  Pick a function
$v \in C^\infty_0(B)$ such that $\Vert v \Vert_{L^2(B)}=1$. Then
\be
\lab{2.115}
L(v, g, \a, B)=L(v, g, 1, B)  + (\a - 1) \frac{n}{2} \ln \left(\int_{B} ( 4 |\nabla v |^2 + R v^2 ) dg \right),
\ee and hence
\be
\lab{2.116}
L(v, g, \a, B) \ge  \lam(g)  + (\a - 1) \frac{n}{2} \ln \left(\int_{B} ( 4 |\nabla v |^2 + R v^2 ) dg \right).
\ee This shows,
\be
\lab{2.117}
L(v, g, \a, B) \ge \lam(g) + (\a - 1) \frac{n}{2} \ln \left( A^{-1} \Vert v \Vert^2_{L^{2n/(n-2)}(B)}
\right),
\ee which implies, via H\"older inequality,
\be
\lab{2.118}
L(v, g, \a, B) \ge -  \frac{n}{2} + (\a - 1) \frac{n}{2} \ln \left( A^{-1} \Vert v \Vert^2_{L^2(B)} /
|B|^{2/n}
\right).
\ee Thus
\be
\lab{2.119}
L(v, g, \a, B) \ge -  \frac{n}{2} - (\a - 1) \frac{n}{2} \ln \left( A \,
|B|^{2/n} \right),
\ee proving part (a) of the proposition. One can also use the fact that
$L(v, g, \a, B) \ge \lambda(g)+(\alpha-1) ( \int_B v^2 \ln v^2 dg -C)$ and
$v^2 \ln v^2 \ge - e^{-1}$ to get the proof.

Now we prove part (b). Notice that in the last paragraph we actually showed that
\[
L(v, g, \a, B) \ge L(v, g, 1, B) - (\a - 1) \frac{n}{2} \ln (A |B|^{2/n}).
\]Hence
\[
\liminf_{\a \to 1^+} \lambda(g, \a, B) \ge \lambda(g, 1, B).
\]

Next we pick, for any given $\epsilon>0$, a function
$v \in C^\infty_0(B)$ such that $\Vert v \Vert_2=1$ and that
\be
\lab{2.120}
\al
\lambda(g, 1, B) &\ge L(v, g, 1, B) -\epsilon\\
&=-\int_B v^2 \ln v^2 + \frac{\alpha}{2} n \ln F(v) + s_n
 + (1-\alpha) \frac{n}{2} \ln F(v) - \epsilon\\
&\ge \lambda(g, \alpha, B) +  (1-\alpha) \frac{n}{2} \ln F(v) -
\epsilon. \eal \ee Since $v$ is fixed, we deduce, after letting
$\alpha \to 0$, that \be \lab{2.121} \lambda(g, 1, B) \ge
\limsup_{\a \to 1^+}\lambda(g, \alpha, B) -\e. \ee Part (b) of the
proposition follows from this when $\e \to 0$. \qed

\begin{proposition}
\lab{prlamtosob}
Let $(M, g)$ be a noncompact manifold  with bounded curvature such that $\lam(g)> -\infty$. If also the scalar curvature
$R \ge 0$,
then there exists a positive constant $A$ depending only on $\lam(g)$ and $n$ such that
\be
\lab{2.122}
\left( \int_{M} v^{2n/(n-2)} dg \right)^{(n-2)/n} \le
A \int_{M} ( 4 |\nabla v |^2 + R v^2 ) dg, \quad \forall  v \in
W^{1, 2}(M, g).
\ee  Moreover, $M$ is $\kappa$ non-collapsed under all scales. i.e. there exists $\kappa>0$ such that
\be
\lab{2.123}
|B(x, r)| \ge \kappa r^n, \qquad r \in (0, \infty)
\ee provided that $R \le 1/r^2$ in $B(x, r)$.
\end{proposition}
\proof
This statement is nothing but the well known equivalence of the Sobolev
inequality and log Sobolev inequality, which is proved via an upper bound for the heat kernel $e^{ (4 \Delta - R) t}$. When $R=0$ one can find a proof in Davies \cite{Da:1}
Chapter 2. When $R \ge 0$, then the
$L^1$ to $L^1$ norm of the heat kernel is less than or equal to
$1$.  The same proof still goes through as written in \cite{Zshu:1} Section 6.2.

Now, we assume $R \le 1/r^2$ in $B(x, r)$. Then
\be
\lab{2.124}
\left( \int_{B(x, r)} v^{2n/(n-2)} dg \right)^{(n-2)/n} \le
A \int_{B(x, r)} ( 4 |\nabla v |^2 + \frac{1}{r^2} v^2 ) dg, \quad \forall  v \in
W^{1, 2}_0(B(x, r), g).
\ee
It is well known that the above Sobolev inequality implies that
 $|B(x, r)| \ge \kappa r^n$ for some $\kappa>0$. See \cite{Ak:1}
 and \cite{Cn:1} e.g.
 Since $x$ and $r$ are arbitrary, $M$ is $\kappa$ noncollapsed under all scales. \qed

As an application of the log Sobolev functional, we next show that
 the Ricci flow on AF manifolds with positive scalar curvature is uniformly $\kappa$ noncollapsed
 for all time.  This result, being different from Perelman's local noncollapsing result
 which holds in finite time, seems to have implications for the issue of longtime convergence.
For example, if the scaled curvature stays bounded, then the
 Gromov-Hausdorf  limit as $t \to \infty$ is still a smooth Riemannian manifold.

\begin{theorem}
\lab{thallnoncol}
Let $(M, g(t))$, $t \in [0, T)$, $T \le \infty$,  be a smooth Ricci flow on AF manifold $M$ of dimension $n \ge
3$. Suppose the scalar curvature $R$ is positive everywhere. Then $(M, g(t))$ is  uniformly $\kappa$
noncollapsed under all scales and for all time.  Moreover,  there exists $A>0$ which depends only
on the initial metric $g(0)$ such that
\be
\lab{alltimeSob}
\left( \int_{M} v^{2n/(n-2)} dg(t) \right)^{(n-2)/n} \le
A \int_{M} ( 4 |\nabla v |^2 + R v^2 ) dg(t), \quad \forall  v \in
W^{1, 2}(M, g(t)), \quad t \in (0, T).
\ee
\end{theorem}
\proof
We just need to prove (\ref{alltimeSob}) since the statement on $\kappa$ noncollapsing follows as mentioned in the previous proposition.

According to Proposition \ref{prAFsob}, $\lam(g(0)) \ge -C > -\infty$. We claim
that $\lam(g(t))$ is monotone nondecreasing in time. Here goes the proof.  Let $t_1, t_2 \in
[0, T)$ and $t_1<t_2$.  For any $\e>0$, there exists a function $\phi \in C^\infty_0(M, g(t_2))$
such that $\Vert \phi \Vert_{L^2(g(t_2))} = 1$ and that
\be
\lab{2.125}
\lam(g(t_2)) \ge L(\phi, g(t_2), 1) - \e.
\ee Now, following Perelman,   let $u=u(x, t)$ be the solution of the conjugate heat equation
with final value $u(x, t_2) = \phi^2(x)$.
Then,  as shown in (\ref{dldt>}) during the proof of Theorem \ref{thnobreather} below,
\be
\lab{2.126}
\frac{d}{dt} L(\sqrt{u(\cdot, t)}, g(t)) \ge 0.
\ee Hence
\be
\lab{2.127}
\lam(g(t_2)) \ge L(\phi, g(t_2), 1) - \e \ge L(\sqrt{u(\cdot, t_1)}, g(t_1), 1) - \e
\ge \lam(g(t_1)) - \e.
\ee This proves the claim and therefore
\be
\lab{2.128}
\lam(g(t)) \ge \lam(g(0)) \ge -C, \quad \forall t>0.
\ee By Proposition \ref{prlamtosob}, we know that (\ref{alltimeSob}) is true.
\qed

\section{Proof of Theorems}

We will prove a number of lemmas first and proceed to prove Theorems \ref{thminimizer} and
\ref{thnobreather}.  During the proof, we will often consider the scaled up manifolds $(M, c_k g, x_k)$ where $c_k \to \infty$ and $x_k$ is a sequence of points in $M$ that may or may
not be fixed.  By the boundedness assumption of the curvature tensor and $\kappa$ noncollapsing condition, we know that this sequence of pointed manifolds sub-converges in $C^\infty_{loc}$ sense, to the Euclidean space with flat metric. This process obviously works for
asymptotically flat manifolds. Notice that the asymptotical
flatness in Corollary \ref{coAF} does not contribute or interfere with this limiting process.
We also do not require that each of the manifold $(M, c_k g, x_k)$ is asymptotically flat
in a uniform way.
The key quantities $\lambda_\infty$ and $\lambda$ only enter the proof through the equations
of the minimizers.

First we show that a minimizer for the functional $L(\cdot, g, \a, B)$ exists when $\a>1$ and
$B$ is a ball.

 \begin{lemma}
\lab{lealphajixiao}
Let $(M, g)$  be a noncompact manifold such that $\lam(g)>-\infty$ and the scalar curvature $R \ge 0$.

(a).  For any
$x_0 \in M$,  and $r_0>0$, write $B=B(x_0, r_0)$. Then for all $\a > 1$,
the infimum of the log Sobolev functionals $L(\cdot, g, \a, B)$ is reached.
Namely, there exists a function $v \in C^\infty_0(B)$ with unit $L^2$ norm such that
\be
\lab{3.99}
L(v, g, \a, B) = \lam(g, \a,  B).
\ee

(b).  The function $v$, called the minimizer, satisfies the equation
\be
\lab{eqlam-a}
\a \frac{n}{2} \frac{ 4 \Delta v - R v}{\int_B ( 4 |\nabla v |^2 + R v^2 ) dg }
+ 2 v \ln v + \b v = 0,
\ee where
\be
\lab{beta=}
\b =  \lam(g, \a,  B)+  \a \frac{n}{2} - \a \frac{n}{2}  \ln \left[ \int_B ( 4 |\nabla v |^2 + R v^2 ) dg \right]
- s_n.
\ee Here $s_n$ is the number given in Definition \ref{deflv}.
\end{lemma}

\proof By Proposition \ref{prlamal}, the  log Sobolev functional is bounded from below.
Hence there exists a sequence of functions $ \{ v_k \} \subset
W^{1, 2}_0(B)$ with unit $L^2$ norm such that
\be
\lab{lvk-lam}
L(v_k, g, \a, B ) \to \lam(g, \a, B) > - \infty, \quad k \to \infty.
\ee So, for all large $k$, we have
\be
\lab{3.100}
-\int_B v^2_k \ln v^2_k dg + \a \frac{n}{2}  \ln \left[ \int_B ( 4 |\nabla v_k |^2 + R v^2_k ) dg \right] + s_n \le
\lam(g, \a, B) + 1.
\ee By the assumption $\lam(g)>\infty$ and Proposition \ref{prlamal},
\be
\lab{3.101}
-\int_B v^2_k \ln v^2_k dg +  \frac{n}{2}  \ln \left[ \int_B ( 4 |\nabla v_k |^2 + R v^2_k ) dg \right]
\ge \lam(g, 1, B) \ge - \lam(g, 1) \le -C  > -\infty.
\ee Substituting this to the previous inequality, we obtain
\be
\lab{fv<}
(\a - 1)  \frac{n}{2}  \ln \left[ \int_B ( 4 |\nabla v_k |^2 + R v^2_k ) dg \right] \le
\lam(g, \a, B) + C - s_n + 1.
\ee By Proposition \ref{prlamtosob}
\be
\lab{3.102}
\al
A^{-1} \left( \int_{B} v^{2n/(n-2)}_k dg \right)^{(n-2)/n} &\le
\int_B ( 4 |\nabla v_k |^2 + R v^2_k ) dg \\
&\le \exp
\left[ (\a-1)^{-1}  \, (\lam(g, \a, B) + C - s_n + 1) \right].
\eal
\ee Pick a number $q \in (2, 2n/(n-2))$. Since the embedding to $L^q(B)$ is compact, we can find
a subsequence, still denoted by $\{ v_k \}$, which converges strongly to a function $v$ in
$L^q(B)$ norm.  By (\ref{fv<}), clearly $v \in W^{1, 2}_0(B)$.

Now we show that $v$ is a minimizer for $L(\cdot, g, \a, B)$. By
Fatou's lemma
\be
\lab{fvk<} \int_B ( 4 |\nabla v |^2 + R v^2 ) dg
\le \lim_{k \to \infty} \int_B ( 4 |\nabla v_k |^2 + R v^2_k ) dg.
\ee According to Theorem 2 in \cite{BL:1},
\be
\lab{3.103}
\int_B v^2 \ln v^2 dg = \lim_{k \to \infty} \int_B v^2_k \ln v^2_k dg +
\lim_{k \to \infty} \int_B (v_k - v)^2 \ln (v_k - v)^2 dg.
\ee Write $B_k = \{ x \, | \, |v_k(x)-v(x)| \le 1 \}$. Then
\be
\lab{3.104}
\int_B (v_k - v)^2 \ln (v_k - v)^2 dg = \int_{B_k} (v_k - v)^2 \ln (v_k - v)^2 dg + \int_{B-B_k} (v_k - v)^2 \ln (v_k - v)^2 dg,
\ee  and therefore
\be
\lab{3.105}
\left| \int_B (v_k - v)^2 \ln (v_k - v)^2 dg \right|
\le \left| \int_{B_k} (v_k - v)^2 \ln (v_k - v)^2 dg  \right| + C_q  \int_{B-B_k} |v_k - v|^q dg.
\ee Applying dominated convergence theorem on the first term of the right hand side, we obtain,
since also $v_k \to v$ in $L^q(B)$ norm, that
\be
\lab{3.106}
\lim_{k \to \infty} \int_B (v_k - v)^2 \ln (v_k - v)^2 dg = 0.
\ee  Consequently
\be
\lab{3.107}
\int_B v^2 \ln v^2 dg = \lim_{k \to \infty} \int_B v^2_k \ln v^2_k dg.
\ee By this and (\ref{fvk<}), we find that
\be
\lab{3.108}
L(v, g, \a, B) \le lim_{k \to \infty} L(v_k, g, \a, B) = \lam(g, \a, B) \le L(v, g, \a, B).
\ee  Hence $v$ is a minimizer.  By the Lagrange multiplier method, there is a constant $\b$ such that
\be
\lab{eqvalph}
\a \frac{n}{2} \frac{ 4 \Delta v - R v}{\int_B ( 4 |\nabla v |^2 + R v^2 ) dg }
+ 2 v \ln v + \b v = 0.
\ee Since $F \equiv \int_B ( 4 |\nabla v |^2 + R v^2 ) dg$ is a finite number, we can multiply it on both sides
of the equation to obtain
\be
\lab{3.109}
\a \frac{n}{2}  4 \Delta v - R v
+ F 2 v \ln v + F \b v = 0.
\ee Since the nonlinear term $v \ln v$ is very mild, it is known that $v \in C^\infty_0(B)$. See \cite{Rot:1}  e.g.

Multiplying (\ref{eqvalph}) by $v$ and integrating, we deduce
\be
\lab{3.110}
-\a \frac{n}{2} + \int_B v^2 \ln v^2 dg + \b =0.
\ee  Since we have proven that $v$ is a minimizer for $L(\cdot, g, \a, B)$, we know that
\be
\lab{3.111}
\lam(g, \a, B) = -\int_B v^2 \ln v^2 dg +\a \frac{n}{2} \ln F + s_n.
\ee Combining the last two identity, we see that
\be
\lab{3.112}
\b =  \lam(g, \a,  B)+  \a \frac{n}{2} - \a \frac{n}{2}  \ln F
- s_n,
\ee  which is just (\ref{beta=}). \qed

The next lemma deals with the case $\a=1$.

\begin{lemma}
\lab{le1jixiao}

Let $(M, g)$  be a noncompact manifold such that $\lam(g)>-\infty$ and that the scalar
curvature $R \ge 0$.

(a).  For any
$x_0 \in M$ and $r_0>0$, let $B=B(x_0, r_0)$. If
$\lam(g, 1, B)<0$, then,
the infimum of the log Sobolev functionals $L(\cdot, g, 1, B)$ is reached.
Namely, there exists a function $v \in C^\infty_0(B)$ with unit $L^2$ norm such that
\be
\lab{3.113}
L(v, g, 1, B) = \lam(g, 1,  B).
\ee

(b).  The function $v$, called the minimizer, satisfies the equation
\be
\lab{eqlam1}
 \frac{n}{2} \frac{ 4 \Delta v - R v}{\int_B ( 4 |\nabla v |^2 + R v^2 ) dg }
+ 2 v \ln v + \b v = 0,
\ee where
\be
\lab{beta1=}
\b =  \lam(g, 1,  B)+   \frac{n}{2} -  \frac{n}{2}  \ln \left[ \int_B ( 4 |\nabla v |^2 + R v^2 ) dg \right]
- s_n.
\ee Here $s_n$ is the number given in Definition \ref{deflv}.
\end{lemma}

\proof  The proof is consisted of a number of steps.

{\it step 1.}  constructing an approximating sequence.

Pick a sequence $\a_k \to 1^{+}$, as $k \to \infty$. Let $v_k$ be a minimizer for $L(\cdot, g, \a_k, B)$,
which exists according to Lemma \ref{lealphajixiao}, and which satisfies
\be
\lab{eqvk}
\a_k \frac{n}{2} \frac{ 4 \Delta v_k - R v_k}{\int_B ( 4 |\nabla v_k |^2 + R v^2_k ) dg }
+ 2 v_k \ln v_k + \b_k v_k = 0,
\ee where
\be
\lab{betak=}
\b_k =  \lam(g, \a_k,  B)+  \a_k \frac{n}{2} - \a_k \frac{n}{2}  \ln \left[ \int_B ( 4 |\nabla v_k |^2 + R v^2_k ) dg \right]
- s_n.
\ee Write
\be
\lab{fkmk}
F_k \equiv \int_B ( 4 |\nabla v_k |^2 + R v^2_k ) dg, \qquad
m_k = \max \{v_k(x) \, | \, x \in B \}.
\ee Since $v_k = 0$ on $\partial B$, we know $\Delta v_k \le 0$ at the maximum point of $v_k$. Hence
(\ref{eqvk}) implies, at the maximum point,
\be
\lab{3.114}
2 v_k \ln v_k \ge -\b_k v_k + \a_k  \frac{n}{2} R v_k F^{-1}_k \ge -\b_k v_k.
\ee By Lemma \ref{prlamal}
\be
\lab{3.115}
\lim_{k \to \infty} \lam(g, \a_k, B)= \lam(g, 1, B)<0.
\ee Therefore, for sufficiently large $k$, we also have $\lam(g, \a_k, B)<0$. This fact and (\ref{betak=})
infer that
\be
\lab{mk>}
m_k = \max v_k \ge e^{- \a_k n/4} F^{\a_k n/4}_k e^{s_n/2}.
\ee

Next we perform the scaling
\be
\lab{3.116}
g_k = m^{4/n}_k g, \quad R_k = m^{-4/n}_k R,  \quad \t v_k = m^{-1}_k v_k.
\ee Notice that $0 \le \t v_k \le 1$ and that
\be
\lab{3.117}
\Vert \t v_k \Vert_{L^2(M, g_k)} = 1.
\ee By (\ref{eqvk}), $\t v_k$ satisfies the equation
\be
\lab{3.118}
\al
\a_k \frac{n}{2} F^{-1}_k m^{4/n}_k& ( 4 \Delta_{g_k} - m^{-4/n}_k R)( m_k \t v_k)
+ 2 m_k \t v_k \ln (m_k \t v_k) \\
&+ ( \lam(g, \a_k,  B)+  \a_k \frac{n}{2} - \a_k \frac{n}{2}  \ln F_k
- s_n) (m_k \t v_k) = 0
\eal
\ee  which becomes, after simplification,
\be
\lab{eqtvk}
\al
\a_k \frac{n}{2}  & ( 4 \Delta_{g_k} -  R_k)  \t v_k
+ ( 2 \t v_k \ln \t v_k  + \lam(g, \a_k,  B) \t v_k +  \a_k \frac{n}{2} \t v_k
- s_n  \t v_k ) \, F_k m^{-4/n}_k \\
&- \a_k \frac{n}{2}  F_k m^{-4/n}_k \ln (F_k m^{-4/(n \a_k)}_k) \,
\t v_k= 0.
\eal
\ee Here $B=B(x_0, r_0, g)$ again.

{\it step 2.}  We prove that for all sequences $\{ \a_k \} \subset (1, 2]$ such that $\a_k \to 1$,
and fixed $r_0$ sufficiently large, there exists a uniform constant $C_0$ such that
\be
\lab{fkjie}
\lim \sup_{ k \to \infty} F_k \le C_0=C_0(r_0).
\ee

Suppose for contradiction that there exists a sequence of
numbers $\{ \a_k \} \subset (1, 2]$ such that $\a_k \to 1$, and that $v_k$ is a minimizer of
$L(\cdot, g, \a_k, B(x_0, r_0))$ but
\be
\lab{fkinfty}
\lim_{ k \to \infty} F_k= \lim_{ k \to \infty} \int_{B(x_0, r_0)} ( 4 |\nabla v_k |^2 + R v^2_k ) dg = \infty.
\ee Then (\ref{mk>}) shows that $m_k \to \infty$ as $k \to \infty$ and that there exists a constant
$C$ such that
\be
\lab{3.119}
F_k m^{-4/(n \a_k)}_k \le C,
\ee and when $k$ is large
\be
\lab{3.120}
F_k m^{-4/n}_k \le F_k m^{-4/(n \a_k)}_k \le C; \qquad a_k \frac{n}{2}  F_k m^{-4/n}_k | \ln (F_k m^{-4/(n \a_k)}_k) | \le C.
\ee Therefore the coefficients of equation (\ref{eqtvk}) are uniformly bounded. Moreover the manifold
$(M, g_k)$ has uniformly bounded geometry since $g_k  = m^{4/n}_k g$ and $m_k \to \infty$.
Now we extend $\t v_k$ to be a function on the whole manifold $M$ by setting $\t v_k =0$
outside of $B(x_0, r_0, g) = B(x_0, m^{2/n}_k r_0, g_k) $. The extended function, still denoted by $\t v_k$,
is a subsolution of the equation in (\ref{eqtvk}); further more $0 \le \t v_k \le 1$ and $\Vert \t v_k \Vert_{L^2(M, g_k)}
=1$.

Let $x_k$ be a maximum point of $\t v_k$ and $r>0$ be a large number. Construct a standard
cut-off function $\phi$ such that $\phi=1$ on $B(x_k, r, g_k)$, $\phi=0$ outside of $B(x_k, 2 r, g_k)$,
$0 \le \phi \le 1$ and $| \nabla_{g_k} \phi | \le C/r$. Since the extended function $\t v_k$ is a sub-solution of (\ref{eqtvk}), we can use $\t v_k \phi^2$ as a test function to conclude, using the
bounds in the previous paragraph, that
\be
\lab{dvk<}
\al
&\int_{B(x_k, r, g_k)} | \nabla_{g_k} \t v_k |^2 d g_k  \\
&\le \frac{C}{r^2} \int_{B(x_k, 2 r, g_k)}  \t v^2_k d g_k
 + C( 1+ |\lam(g, \a_k,  B )| )
F_k m^{-4/n}_k
 \int_{B(x_k, 2 r, g_k)}  \t v^2_k d g_k\\
 & \le \frac{C}{r^2}
 + C( 1+ |\lam(g, \a_k,  B ) | )
F_k m^{-4/n}_k.
\eal
\ee Here $B=B(x_0, r_0, g)$ again.

We consider 2 cases.
\medskip

Case 1.  A subsequence of $\{ F_k m^{-4/n}_k \}$,  denoted by the same symbol, converges to $0$.

Let $x_k$ be a maximum point of $v_k$. Since $m_k \to \infty$ and $g_k = m^{4/n}_k g$,
we know that a subsequence of the pointed manifolds $\{(M, g_k, x_k)\}$, converges in $C^\infty_{loc}$
topology, to the pointed Euclidean space $({\R^n}, 0)$.  This is due to the Cheeger-Gromov compactness theorem.
By the bound (\ref{dvk<}) and the fact $R_k \to 0$, $\lam(g, \a_k,  B ) \to \lam(g, 1, B)$, we know that
a subsequence of $\t v_k$ converges pointwise, modulo composition with diffeomorphisms, to a
function $v_\infty$ on $\R^n$, which is a sub-solution of the Laplacian. Furthermore
$\Vert v_\infty \Vert_{L^2(\R^n)} \le 1$ and  $v_\infty(0)=1$. By (\ref{dvk<}) again
\be
\lab{3.121}
\int_{B(0, r)} | \nabla  v_\infty |^2 dx \le \frac{C}{r^2}.
\ee Here all expressions are in the Euclidean setting.  Letting $r \to \infty$, we see that $\nabla v_\infty
=0$ and therefore $v_\infty \equiv 1$. But this is impossible since $\Vert v_\infty \Vert_{L^2(\R^n)} \le 1$.

\medskip

Case 2.  $\{ F_k m^{-4/n}_k \}$ is bounded away from $0$.

Then we can find
 a subsequence of $\{ F_k m^{-4/n}_k \}$,  denoted by the same symbol, which converges to a number $A>0$. As in the previous paragraph, $\{(M, g_k, x_k)\}$, converges in $C^\infty_{loc}$
topology, to the pointed Euclidean space $({\R^n}, 0)$.
Also a subsequence of the extended function $\t v_k$ converges pointwise, modulo composition with diffeomorphisms, to a
function $v_\infty$ on $\R^n$.   Furthermore
$\Vert v_\infty \Vert_{L^2(\R^n)} \le 1$ , $v_\infty(0)=1$ and, in the weak sense,
\be
\lab{eqabv}
\frac{n}{2}   4 \Delta v_\infty
+ A ( 2 v_\infty \ln  v_\infty  + \lam(g, 1,  B) v_\infty +  \frac{n}{2} v_\infty
- s_n  v_\infty)  \\
-  (\frac{n}{2}  A \ln A) \,
 v_\infty \ge 0.
\ee Dividing both sides by $A$ and recalling  from Definition \ref{deflv} that $s_n =  -\frac{n}{2} \ln (2 \pi n) - \frac{n}{2}$, we obtain
\be
\lab{3.122}
 \lam(g, 1,  B) v_\infty \ge - \frac{n}{2A}   4 \Delta v_\infty  - 2 v_\infty \ln  v_\infty -  n v_\infty
 +\frac{n}{2}  \ln (2 \pi n A) \,
 v_\infty.
\ee

We multiply the last inequality by $v_\infty$.  By Moser's iteration, it is easy to prove that  $v_\infty$ has
Gaussian decay near infinity. See \cite{Rot:1} or Lemma 2.3 in \cite{Z:2} e.g.  Therefore, we can carry out integration by parts to deduce
\be
\lab{3.123}
\al
 \lam(g, 1,  B) &  \Vert v_\infty \Vert^2_{L^2(\R^n)} \ge \int_{\R^n}  \left( \frac{n}{2A}   4 |\nabla v_\infty|^2  -  v^2_\infty \ln  v^2_\infty -  n v^2_\infty
 -\frac{n}{2}  \ln (2 \pi n /A) \,
 v^2_\infty \right) dx\\
 &=\int_{\R^n}  \left( s   4 |\nabla v_\infty|^2  -  v^2_\infty \ln  v^2_\infty
 -\frac{n}{2}  \ln (4 \pi s) \,
 v^2_\infty -  n v^2_\infty \right) dx,
 \eal
\ee where $s= \frac{n}{2A}$. Write $\hat v =  \frac{v_\infty}{\Vert v_\infty \Vert_{L^2(\R^n)}}$. Then,
by $\Vert  v_\infty \Vert^2_{L^2(\R^n)} \le 1$, we have
\be
\lab{3.124}
\al
 \lam(g, 1,  B)  \Vert v_\infty \Vert^2_{L^2(\R^n)}
&\ge \Vert  v_\infty \Vert^2_{L^2(\R^n)}  \int_{\R^n}  \left( s   4 |\nabla \hat v|^2  -  \hat v^2 \ln  \hat v^2_\infty
 -\frac{n}{2}  \ln (4 \pi s) \,
 \hat v^2_\infty -  n \hat v^2_\infty \right) dx  \\
 &\qquad - \Vert v_\infty \Vert^2_{L^2(\R^n)}  \ln \Vert  v_\infty \Vert^2_{L^2(\R^n)}  \ge 0.
 \eal
\ee
Here we just used the fact that the best constant for the log Sobolev inequality for functions with unit $L^2$ norms in $\R^n$ is $0$.  This is a contradiction with the assumption that $\lam(g, 1,  B)<0$.
This proves (\ref{fkjie}), i.e.
\be
\lab{3.125}
F_k = \int_{B(x_0, r_0)} ( 4 |\nabla v_k |^2 + R v^2_k ) dg \le C_0.
\ee

{\it step 3.}  We prove $v_k$ converges to a minimizer of $L(\cdot, g, 1, B)$.
\medskip

By (\ref{eqvk}), we know that $v_k$ satisfies
\be
\lab{3.126}
 \frac{n}{2} ( 4 \Delta v_k - R v_k)
+ \a_k^{-1} F_k 2 v_k \ln v_k + \a_k^{-1} F_k \b_k v_k = 0,
\ee where
\be
\lab{3.127}
\b_k =  \lam(g, \a_k,  B)+  \a_k \frac{n}{2} - \a_k \frac{n}{2}  \ln F_k
- s_n.
\ee Since, by Step 2, $F_k$ is uniformly bounded, we know that the coefficients in the above equation
are uniformly bounded.  Since the nonlinear term $v_k \ln v_k$
is only mildly nonlinear, it is easy to prove that $\Vert v_k \Vert_{L^\infty}$ is
also uniformly bounded. See Lemma 2.1 in \cite{Z:2} e.g. Now, since  the ball $B$ is bounded,a routine argument shows that a subsequence of $v_k$
converges to a minimizer $v$ of $L(\cdot, g, 1, B)$.  Using the same argument near the end of the
proof of Lemma \ref{lealphajixiao}, we see that $v$ satisfies equation (\ref{eqlam1}).
This proves the lemma .
\qed

 \medskip

 The next lemma shows that the minimizers of $L(\cdot, g, 1, B)$  are uniformly bounded
 even if the radius of $B$ tends to $\infty$.

 \begin{lemma}
\lab{leyizhijie}
Under the same assumption as in Lemma \ref{le1jixiao},
let $v$ be a minimizer for $L(\cdot, g, 1, B)$, where $B=B(x_0, r_0)$.
Then the quantity
\be
\lab{3.128}
F = \int_B ( 4 |\nabla v|^2 + R v^2) dg
\ee is uniformly bounded for all large $r_0$. Furthermore $ \Vert v \Vert_{L^\infty(B)}$ is uniformly
bounded for all large $r_0$.
\end{lemma}

\proof The idea of the proof is similar to that for the previous lemma.
Suppose for contradiction that there exists a sequence of radii $\{ r_{0k} \}$ and that $v_k$ is a minimizer of
$L(\cdot, g, 1, B(x_0, r_{0k}))$ but
\be
\lab{fkinfty2}
\lim_{ k \to \infty} F_k= \lim_{ k \to \infty} \int_{B(x_0, r_{0k})} ( 4 |\nabla v_k |^2 + R v^2_k ) dg = \infty.
\ee From the previous lemma $v_k$ satisfies
\be
\lab{eq1vk}
 \frac{n}{2} \frac{ 4 \Delta v_k - R v_k}{F_k }
+ 2 v_k \ln v_k + \b_k v_k = 0,
\ee where
\be
\lab{beta1k=}
\b_k =  \lam(g, 1,  B_k)+   \frac{n}{2} -  \frac{n}{2}  \ln F_k
- s_n.
\ee Here and later $B_k =  B(x_0, r_{0k}) = B(x_0, r_{0k}, g)$.
 Since $v_k = 0$ on $\partial B_k$, we know $\Delta v_k \le 0$ at the maximum point of $v_k$. Hence
(\ref{eq1vk}) implies, at the maximum point,
\be
\lab{3.129}
2 v_k \ln v_k \ge -\b_k v_k +   \frac{n}{2} R v_k F^{-1}_k \ge -\b_k v_k.
\ee Since by definition
\be
\lab{3.130}
\lim_{k \to \infty} \lam(g, 1, B_k)= \lam(g, 1, M)<0,
\ee for sufficiently large $k$, we also have $\lam(g, 1, B_k)<0$. This fact and (\ref{beta1k=})
infer that
\be
\lab{1mk>}
m_k \equiv \max v_k \ge e^{-  n/4} F^{ n/4}_k e^{s_n/2}.
\ee

Next we do the scaling
\be
\lab{3.131}
g_k = m^{4/n}_k g, \quad R_k = m^{-4/n}_k R,  \quad \t v_k = m^{-1}_k v_k.
\ee Notice that $0 \le \t v_k \le 1$ and that
\be
\lab{3.132}
\Vert \t v_k \Vert_{L^2(M, g_k)} = 1.
\ee By (\ref{eq1vk}), $\t v_k$ satisfies the equation
\be
\lab{3.133}
\al
 \frac{n}{2} F^{-1}_k m^{4/n}_k& ( 4 \Delta_{g_k} - m^{-4/n}_k R)( m_k \t v_k)
+ 2 m_k \t v_k \ln (m_k \t v_k) \\
&+ ( \lam(g, 1,  B_k)+   \frac{n}{2} -  \frac{n}{2}  \ln F_k
- s_n) (m_k \t v_k) = 0
\eal
\ee  which becomes, after simplification,
\be
\lab{eq1tvk}
\al
 \frac{n}{2}  & ( 4 \Delta_{g_k} -  R_k)  \t v_k
+ ( 2 \t v_k \ln \t v_k  + \lam(g, 1,  B_k) \t v_k +   \frac{n}{2} \t v_k
- s_n  \t v_k ) \, F_k m^{-4/n}_k \\
&-  \frac{n}{2}  F_k m^{-4/n}_k \ln (F_k m^{-4/n}_k) \,
\t v_k= 0.
\eal
\ee

Since $F_k \to \infty$ by assumption,
 (\ref{1mk>}) shows that $m_k \to \infty$ as $k \to \infty$ and that there exists a constant
$C$ such that
\be
\lab{3.134}
F_k m^{-4/n}_k \le C,
\ee Therefore the coefficients of equation (\ref{eq1tvk}) are uniformly bounded. Moreover the manifold
$(M, g_k)$ has uniformly bounded geometry since $g_k  = m^{4/n}_k g$ and $m_k \to \infty$.
Now we extend $\t v_k$ to be a function on the whole manifold $M$ by setting $\t v_k =0$
outside of $B_k=B(x_0, r_{0k}, g) = B(x_0, m^{2/n}_k r_{0k}, g_k) $. The extended function, still denoted by $\t v_k$,
is a subsolution of (\ref{eq1tvk}); further more $0 \le \t v_k \le 1$ and $\Vert \t v_k \Vert_{L^2(M, g_k)}
=1$.

Let $x_k$ be a maximum point of $\t v_k$ and $r>0$ be a large number. Construct a standard
cut-off function $\phi$ such that $\phi=1$ on $B(x_k, r, g_k)$, $\phi=0$ outside of $B(x_k, 2 r, g_k)$,
$0 \le \phi \le 1$ and $| \nabla_{g_k} \phi | \le C/r$. Since the extended function $\t v_k$ is a sub-solution of (\ref{eq1tvk}), we can use $\t v_k \phi^2$ as a test function to conclude, using the
bounds in the previous paragraph, that
\be
\lab{d1vk<}
\al
&\int_{B(x_k, r, g_k)} | \nabla_{g_k} \t v_k |^2 d g_k  \\
&\le \frac{C}{r^2} \int_{B(x_k, 2 r, g_k)}  \t v^2_k d g_k
 + C( 1+ |\lam(g, 1,  B_k )| )
F_k m^{-4/n}_k
 \int_{B(x_k, 2 r, g_k)}  \t v^2_k d g_k\\
 & \le \frac{C}{r^2}
 + C( 1+ |\lam(g, 1,  B_k )|  )
F_k m^{-4/n}_k.
\eal
\ee

We consider 2 cases.
\medskip

Case 1.  A subsequence of $\{ F_k m^{-4/n}_k \}$,  denoted by the same symbol, converges to $0$.

Let $x_k$ be a maximum point of $v_k$ again. Since $m_k \to \infty$ and $g_k = m^{4/n}_k g$,
by Cheeger-Gromov compactness theorem,
we know that a subsequence of the pointed manifolds $\{(M, g_k, x_k)\}$, converges in $C^\infty_{loc}$
topology, to the pointed Euclidean space $({\R^n}, 0)$.
By the bound (\ref{d1vk<}) and the fact $R_k \to 0$, $\lam(g, 1,  B_k ) \to \lam(g, 1, M)$, we know that
a subsequence of $\t v_k$ converges pointwise, modulo composition with diffeomorphisms, to a
function $v_\infty$ on $\R^n$, which is a sub-solution of the Laplacian. Furthermore
$\Vert v_\infty \Vert_{L^2(\R^n)} \le 1$ and  $v_\infty(0)=1$. By (\ref{d1vk<}) again
\be
\lab{3.135}
\int_{B(0, r)} | \nabla  v_\infty |^2 dx \le \frac{C}{r^2}.
\ee Here all expressions are in the Euclidean setting.  Letting $r \to \infty$, we see that $\nabla v_\infty
=0$ and therefore $v_\infty \equiv 1$. But this is impossible since $\Vert v_\infty \Vert_{L^2(\R^n)} \le 1$.

\medskip

Case 2.  $\{ F_k m^{-4/n}_k \}$ is bounded away from $0$.

Then we can find
 a subsequence of $\{ F_k m^{-4/n}_k \}$,  denoted by the same symbol, which converges to a number $A>0$. As in the previous paragraph, $\{(M, g_k, x_k)\}$, converges in $C^\infty_{loc}$
topology, to the pointed Euclidean space $({\R^n}, 0)$.
Also a subsequence of the extended function $\t v_k$ converges pointwise, modulo composition with diffeomorphisms, to a
function $v_\infty$ on $\R^n$.   Furthermore
$\Vert v_\infty \Vert_{L^2(\R^n)} \le 1$ , $v_\infty(0)=1$ and, in the weak sense,
\be
\lab{eq1abv}
\frac{n}{2}   4 \Delta v_\infty
+ A ( 2 v_\infty \ln  v_\infty  + \lam(g, 1,  M) v_\infty +  \frac{n}{2} v_\infty
- s_n  v_\infty)  \\
-  (\frac{n}{2}  A \ln A) \,
 v_\infty \ge 0.
\ee Dividing both sides by $A$ and recalling  from Definition \ref{deflv} that $s_n =  -\frac{n}{2} \ln (2 \pi n) - \frac{n}{2}$, we obtain
\be
\lab{3.136}
 \lam(g, 1,  M) v_\infty \ge - \frac{n}{2A}   4 \Delta v_\infty  - 2 v_\infty \ln  v_\infty -  n v_\infty
 +\frac{n}{2}  \ln (2 \pi n A) \,
 v_\infty.
\ee

We multiply the last inequality by $v_\infty$.  By Moser's iteration, it is easy to prove,
as in Lemma 2.3 in \cite{Z:2}, $v_\infty$ has
Gaussian decay near infinity. Therefore, we can carry out integration by parts to deduce
\be
\lab{3.137}
\al
 \lam(g, 1,  M) &  \Vert v_\infty \Vert^2_{L^2(\R^n)} \ge \int_{\R^n}  \left( \frac{n}{2A}   4 |\nabla v_\infty|^2  -  v^2_\infty \ln  v^2_\infty -  n v^2_\infty
 -\frac{n}{2}  \ln (2 \pi n /A) \,
 v^2_\infty \right) dx\\
 &=\int_{\R^n}  \left( s   4 |\nabla v_\infty|^2  -  v^2_\infty \ln  v^2_\infty
 -\frac{n}{2}  \ln (4 \pi s) \,
 v^2_\infty -  n v^2_\infty \right) dx,
 \eal
\ee where $s= \frac{n}{2A}$. Write $\hat v =  \frac{v_\infty}{\Vert v_\infty \Vert_{L^2(\R^n)}}$. Then,
by $\Vert  v_\infty \Vert^2_{L^2(\R^n)} \le 1$, we have
\be
\lab{3.138}
\al
 \lam(g, 1,  M)  \Vert v_\infty \Vert^2_{L^2(\R^n)}
&\ge \Vert  v_\infty \Vert^2_{L^2(\R^n)}  \int_{\R^n}  \left( s   4 |\nabla \hat v|^2  -  \hat v^2 \ln  \hat v^2
 -\frac{n}{2}  \ln (4 \pi s) \,
 \hat v^2 -  n \hat v^2 \right) dx  \\
 &\qquad - \Vert v_\infty \Vert^2_{L^2(\R^n)}  \ln \Vert  v_\infty \Vert^2_{L^2(\R^n)}  \ge 0.
 \eal
\ee
Here we just used the fact that the best constant for the log Sobolev inequality for functions with unit $L^2$ norms in $\R^n$ is $0$.  This is a contradiction with the assumption that $\lam(g, 1,  M)<0$.
This proves that $F_k$ is uniformly bounded.

The uniform boundedness of $v_k$ comes from the following arguments.
By (\ref{eq1vk}), we know that $v_k$ satisfies
\be
\lab{3.139}
 \frac{n}{2} ( 4 \Delta v_k - R v_k)
+  F_k 2 v_k \ln v_k +  F_k \b_k v_k = 0,
\ee where
\be
\lab{3.140}
\b_k =  \lam(g, 1,  B_k)+   \frac{n}{2} -  \frac{n}{2}  \ln F_k
- s_n.
\ee Since $F_k$ is uniformly bounded, we know that the coefficients in the above equation
are uniformly bounded. As explained at the end of the proof of Lemma \ref{le1jixiao},  it is easy to
show that $\Vert v_k \Vert_{L^\infty}$ is
also uniformly bounded.  This proves the lemma.
\qed

\medskip

Now we are ready to give

 \proof {\bf of Theorem \ref{thminimizer}.}

We will use the minimizers $v_k$ on balls of radius $r_k$ to construct a minimizer on the
whole manifold. The core argument is to show that $v_k$ has a non vanishing limit.

{\it Step 1.}

 Pick $r_k \to \infty$ and let $v_k$ be a minimizer for $L(\cdot, g, 1, B(x_0, r_k) )$ whose
 infimum is $\lam_k$. Then
 \be
 \lab{lamk=1}
 \al
 \lam_k&=L(v_k, g, 1, B(x_0, r_k) )
\\
&= - \int_{B(x_0, r_k)} v^2_k \ln v^2_k dg + \frac{n}{2} \ln
\left( \int_{B(x_0,
 r_k)} (4 |\nabla v_k|^2 + R v^2_k ) dg \right) + s_n.
 \eal
 \ee
According to the previous 2 lemmas, $v_k$ exists and  is uniformly
bounded. By standard elliptic theory,  a subsequence of $\{ v_k
\}$, still denoted by the same symbol, converges in
$C^\infty_{loc}$ sense, to a limit function $v_\infty \in
C^\infty(M)$.  In this step, we prove that $v_\infty$ is not $0$.
We will use P. L. Lion's concentrated compactness method at infinity.
But a new twist occurs. That is, even though $\lam_k$ is bounded, the
components on the right hand side of (\ref{lamk=1}) may not be bounded
from below uniformly.

Suppose for contradiction that $v_\infty=0$.  Then $v_k \to 0$
a.e. as $k \to \infty$.  Then  there exists a sequence of positive
integers $\{ i_k \}$ and a subsequence of $\{ v_k \}$,  denoted by the same
symbol,
such that $i_k \to \infty$ as $k \to \infty$
and that
\be
\lab{vkBto0} \int_{B(x_0,  2^{2 i_k})}  v^2_k dg \to
0, \qquad k \to \infty.
\ee For any positive integer $i$ we
introduce the following notations
\be
\lab{3.141}
\al
\Omega_i &=B(x_0, 2^i) - B(x_0, 2^{i-1}),\\
F(v_k) & = \int_M ( 4 |\nabla v_k |^2 + R v_k^2 ) dg, \quad
N(v_k) =  \int_M v^2_k \ln v^2_k dg.
\eal
\ee Here $v_k$ is considered $0$ outside of the ball $B(x_0, r_k)$.

By $\lam \equiv \lam(g)= \lam(g, 1)>-\infty$ in assumption (a) of
the theorem and Proposition \ref{prlamtosob}, there exists a positive constant $A$ such that
\be
\lab{3.142}
\left( \int_{B(x_0, r_k)} v^{2n/(n-2)}_k dg \right)^{(n-2)/n} \le
A F(v_k).
\ee

Hence \be \lab{Fv<eNv} \al &\left( \Sigma^{2 i_k}_{i=i_k}
\int_{\Omega_i} v^{2n/(n-2)}_k dg \right)^{(n-2)/n} e^{-N(v_k)
2/n} \le \left( \int_{B(x_0, r_k)}
v^{2n/(n-2)}_k dg \right)^{(n-2)/n} e^{-N(v_k) 2/n}  \\
&\le C F(v_k) e^{-N(v_k) 2/n} = C e^{(\lam_k -s_n) 2/n} \le C,
 \eal
\ee where we also used (\ref{lamk=1}) and the fact that $\lam_k$ is
uniformly bounded. Thus, there exists an integer $j_k \in [i_k, 2
i_k]$ such that
\be \lab{vkojk<} \left(  \int_{\Omega_{j_k}}
v^{2n/(n-2)}_k dg \right)^{(n-2)/n} \le C i_k^{-(n-2)/n}
e^{N(v_k) 2/n} \ee

By partition of unity, we can choose a sequence of cut-off
functions $\phi_k$, $\eta_k$ on $M$ such that $ \phi_k = 1 $ on
$B(x_0, 2^{j_k -1})$,  \, $supp \, \phi_k \subset B(x_0,
2^{j_k})$; $ \eta_k = 1 $ on $M-B(x_0, 2^{j_k})$,  \, $supp \,
\eta_k \subset M- B(x_0,  2^{j_k-1})$;
 $| \nabla \phi_k | + | \nabla \phi_k | \le C/{2^{j_k}}$;
$\phi^2_k + \eta^2_k=1$.  We introduce the notations
\be
\lab{3.143}
a_k \equiv \Vert v_k \phi_k \Vert^2_{L^2}, \quad b_k \equiv \Vert
v_k \eta_k \Vert^2_{L^2};
\ee
\be
\lab{3.144}
A_k \equiv \exp( \frac{2}{n} N(v_k \phi_k)), \quad B_k \equiv
\exp( \frac{2}{n} N(v_k \eta_k)).
\ee By (\ref{vkBto0}), we know that
\be
 \lab{akbkto01} a_k \to 0, \quad b_k \to 1, \quad \text{as}
\quad k \to \infty.
\ee

Now we will split the terms in the log
Sobolev functional into terms involving $v_k \phi_k$ and $v_k
\eta_k$. By direct computation
\be
\lab{fenF}
\al
&\int (4 | \nabla v_k |^2 + R v^2_k ) dg \\
&=\int (4 | \nabla (v_k \phi_k) |^2 + R (v_k \phi_k)^2 ) dg +
\int (4 | \nabla (v_k \eta_k) |^2 + R (v_k \eta_k)^2 ) dg\\
&\qquad - 4 \int (  | \nabla \phi_k|^2 + | \nabla \eta_k|^2 )  v^2_k dg,
\eal
\ee where we have used the identity
\be
\lab{3.145}
0=\Delta (\phi^2_k+\eta^2_k) = 2  | \nabla \phi_k|^2  + 2 \phi_k \Delta \phi_k + 2 | \nabla \eta_k|^2
+ 2 \eta_k \Delta \eta_k.
\ee

Suppose Condition (b) on volume of geodesic balls holds, namely
$|B(x_0, r)| \le C r^n$. Using H\"older's inequality we deduce
\be
\lab{3.146}
\al
4 \int (  | \nabla \phi_k|^2 + | \nabla \eta_k|^2 )  v^2_k dg &\le
C 2^{-2 j_k} \int_{\Omega_{j_k}} v^2_k dg
\le C 2^{-2 j_k} |\Omega_{j_k}|^{2/n} \left(  \int_{\Omega_{j_k}} v^{2n/(n-2)}_k dg \right)^{(n-2)/n}\\
&\le C  \left(  \int_{\Omega_{j_k}} v^{2n/(n-2)}_k dg \right)^{(n-2)/n}.
\eal.
\ee  Using (\ref{vkojk<}), we know that
\be
\lab{dphi+deta} 4 \int (  | \nabla \phi_k|^2 + | \nabla
\eta_k|^2 ) v^2_k dg = o(1)
 e^{N(v_k) 2/n}.
\ee Here $o(1)$ is a quantity that goes to $0$ when $k \to
\infty$. This and (\ref{fenF}) imply
\be \lab{splitFvk} F(v_k) =
F(v_k \phi_k)  + F(v_k \eta_k) - o(1)
 e^{N(v_k) 2/n}.
\ee

Now, suppose Condition (b) on  the scalar curvature holds, namely
$R(x) \ge \frac{c}{1+d^2(x_0, x)}$. Then
\be
\lab{3.147}
4 \int (  | \nabla \phi_k|^2 + | \nabla \eta_k|^2 )  v^2_k dg \le
C 2^{-2 j_k} \int_{\Omega_{j_k}} v^2_k dg \le C
\int_{\Omega_{j_k}} R v^2_k dg.
\ee By the second line of (\ref{Fv<eNv}), we have
\be
\lab{3.148}
\Sigma^{2 i_k}_{i=i_k} \int_{\Omega_i} (4 |\nabla v_k |^2 + R v^2_k)
dg \le C e^{ 2 N(v_k)/n}.
\ee Therefore one can also find a $j_k \in [i_k, 2 i_k]$ such that
(\ref{dphi+deta}) and (\ref{splitFvk}) hold.

Next,
observe that
\be
\lab{3.149}
\al
\int v^2_k &\ln v^2_k dg -
\int (v_k \phi_k)^2 \ln (v_k \phi_k)^2 dg - \int (v_k \eta_k)^2 \ln (v_k \eta_k)^2dg\\
&=\int (v_k \phi_k)^2 \left[ \ln ( (v_k \phi_k)^2 +
(v_k \eta_k)^2 ) -  \ln (v_k \phi_k)^2 \right] dg  \\
&\qquad + \int (v_k \eta_k)^2 \left[ \ln ( (v_k \phi_k)^2  +
(v_k \eta_k)^2 ) -  \ln (v_k \eta_k)^2 \right] dg \\
&\le C \int v^4_k \phi^2_k \eta^2_k dg \le C \int_{\Omega_{j_k}} v^2_k \, dg. \eal
\ee Here we just used the uniform boundedness of $v_k$, proven in Lemma \ref{leyizhijie}.
This means
\be
\lab{splitNvk}
N(v_k) = N(v_k \phi_k) + N(v_k \eta_k) + o(1).
\ee

Recall that $v_k$ is a minimizer for the log Sobolev functional. By (\ref{lamk=1}),
\be
\lab{Fv/eN}
e^{\frac{2}{n} (\lam_k -s_n)} = \frac{F(v_k)}{\exp (\frac{2}{n}
N(v_k)) }.
\ee By (\ref{splitFvk}) and (\ref{splitNvk}), this implies
\be
\lab{3.150}
\al
e^{\frac{2}{n} (\lam_k-s_n)}&=\frac{F(v_k \phi_k) + F(v_k \eta_k)
+ o(1) \exp( \frac{2}{n} N(v_k))}{\exp( \frac{2}{n} N(v_k))} \\
&=\frac{F(v_k \phi_k) + F(v_k \eta_k)
 }{\exp( \frac{2}{n} N(v_k \phi_k)) \, \exp( \frac{2}{n}
N(v_k \eta_k)) \, e^{o(1)}} +o(1).
\eal
\ee  On the other hand, by definition of $\lam_k$, we have
\be
\lab{3.151}
F(v_k \phi_k) \ge e^{\frac{2}{n} (\lam_k-s_n)} \Vert v_k \phi_k
\Vert^2_{L^2} \exp \left( -\frac{2}{n} \ln \Vert v_k \phi_k
\Vert^2_{L^2} \right) \exp\left(\frac{2}{n} N(v_k \phi_k) / \Vert
v_k \phi_k \Vert^2_{L^2}\right).
\ee Since the support of $\eta_k$ is outside of the ball $B(x_0, 2^{j_k-1})$, by definition of
$\lam_\infty \equiv \lam_\infty(g)$ in Definition \ref{deflv}, we
know
\be
\lab{3.152}
F(v_k \eta_k) \ge e^{\frac{2}{n} (\lam_\infty-s_n + o(1))} \Vert v_k \eta_k
\Vert^2_{L^2} \exp \left( -\frac{2}{n} \ln \Vert v_k \eta_k
\Vert^2_{L^2} \right) \exp\left(\frac{2}{n} N(v_k \eta_k) / \Vert
v_k \eta_k \Vert^2_{L^2}\right).
\ee Write $\lam = \lam(g, 1)$. Combining the last three expressions, we deduce, since $\lam_k = \lam+ o(1)$ that
\be
\lab{3.153}
1 \ge \frac{a^{-2/n}_k a_k A^{1/a_k}_k + b^{-2/n}_k b_k B^{1/b_k}_k
e^{(\lam_\infty-\lam) 2/n + o(1)} }{A_k B_k e^{o(1)}}  + o(1),
\ee  where
\be
\lab{3.154}
a_k \equiv \Vert v_k \phi_k \Vert^2_{L^2}, \quad b_k \equiv \Vert v_k \eta_k \Vert^2_{L^2};
\ee
\be
\lab{3.155}
A_k \equiv \exp( \frac{2}{n} N(v_k \phi_k)), \quad B_k \equiv \exp( \frac{2}{n} N(v_k \eta_k)).
\ee  Therefore
\be
\lab{3.156}
\min\{ a^{-2/n}_k, \,b^{-2/n}_k \} \,
\frac{ a_k A^{1/a_k}_k +  b_k B^{1/b_k}_k e^{(\lam_\infty-\lam) 2/n + o(1)} }{A_k B_k e^{o(1)}}  + o(1) \le 1,
\ee Since $a_k$ and $b_k$ are positive numbers in the interval $(0, 1)$, this shows
\be
\lab{3.157}
\ln (a_k A^{1/a_k}_k +  b_k B^{1/b_k}_k e^{(\lam_\infty-\lam) 2/n + o(1)}) \le \ln (A_k B_k) + o(1).
\ee
 Notice that $a_k + b_k =1$. By concavity of $\ln$ function we obtain
\be
\lab{3.158}
b_k (\lam_\infty-\lam) 2/n + o(1) \le o(1).
\ee
Letting $k \to \infty$ and using the fact that $b_k \to 1$ (from (\ref{akbkto01}) ),  we arrive at
\be
\lab{3.159}
0< \lam_\infty-\lam \le 0.
\ee  This is a contradiction  which proves that $v_\infty$ is not identically zero.

\medskip
{\it Step 2.} We prove
$\Vert v_\infty \Vert_{L^2(M)} =1$.

 This is done by adopting a method by Dolbeault and Esteban \cite{DE:1}, which is in the spirit of
P. L. Lions' concentrated compactness.

Suppose for contradiction that $\Vert v_\infty \Vert_{L^\infty(M)} = \delta <1$.
 Then for all large integer $k$, there exists $l_k>0$ such that $l_k \to
\infty$ when $k \to \infty$ and
\be
\lab{3.160}
\int_{B(x_0, l_k)} v^2_\infty dg = \delta - \frac{1}k, \qquad
\int_{B(x_0, 4 l_k)-B(x_0, l_k)} v^2_\infty dg \le \frac{1}k.
\ee  Fixing this $k$ for the moment, by $C^\infty_{loc}$ convergence of $v_k$ to $v_\infty$ and the fact that
the $L^2$ norm of $v_k$ is $1$, we can find
a subsequence $\{ n_k \}$ of positive integers so that
\be
\lab{3.161}
\delta - \frac{2}{k} \le \int_{B(x_0, l_k)} v^2_{n_k} dg \le \delta - \frac{1}{2 k}, \quad
\int_{B(x_0, 4 l_k)-B(x_0, l_k)} v^2_{n_k} dg \le \frac{2}k,
\ee and that
\be
\lab{3.162}
1-\delta - \frac{2}k \le \int_{M-B(x_0, 4 l_k)} v^2_{n_k} dg \le 1-\delta + \frac{2}k.
\ee Renaming $n_k$ as $k$, we have found a subsequence of $\{ v_k \}$, which is still denoted by
$\{ v_k \}$, such that
\be
\lab{v2k=delta}
\al
&\lim_{k \to \infty} \int_{B(x_0, l_k)} v^2_k dg = \delta, \qquad
\lim_{k \to \infty} \int_{B(x_0, 4 l_k)-B(x_0, l_k )} v^2_k dg=0,\\
&\lim_{k \to \infty} \int_{M-B(x_0, 4 l_k)} v^2_k dg = 1-\delta .
\eal
\ee By partition of unity, we can choose a sequence of cut-off functions $\phi_k$, $\eta_k$ on
$(M, x_0, g )$ such that
$
\phi_k = 1
$ on $B(x_0, l_k)$,  \, $supp \, \phi_k \subset B(x_0, 2 l_k)$;
$
\eta_k = 1
$ on $M-B(x_0, 2 l_k)$,  \,  $supp \, \eta_k \subset M- B(x_0,  l_k)$;
 $| \nabla \phi_k | + | \nabla \phi_k | \le C/{l_k}$;
$\phi^2_k + \eta^2_k=1$. Using (\ref{v2k=delta}), we know that
\be
\lab{splitv2}
\lim_{k \to \infty} \int_{B(x_0, l_k)} (v_k \phi_k)^2 dg= \delta, \qquad
\lim_{k \to \infty} \int_{M-B(x_0, 4 l_k)} (v_k \eta_k)^2 dg = 1-\delta .
\ee

Next we will again split the terms in the log Sobolev functional into terms involving $v_k \phi_k$ and
$v_k \eta_k$.
Since $| \nabla \phi_k | + | \nabla \eta_k | \to 0$ when $k \to \infty$, it is easy to see that
\be
\lab{splitF}
\al
&\int (4 | \nabla v_k |^2 + R v^2_k ) dg \\
&=\int (4 | \nabla (v_k \phi_k) |^2 + R (v_k \phi_k)^2 ) dg + \int (4 | \nabla (v_k \eta_k) |^2 + R (v_k
\eta_k)^2 ) dg + o(1).
\eal
\ee Here $o(1)$ is a quantity that goes to $0$ when $k \to \infty$.
As in Step 1,
\be
\lab{3.163}
\al
\int v^2_k &\ln v^2_k dg -
\int (v_k \phi_k)^2 \ln (v_k \phi_k)^2 dg- \int (v_k \eta_k)^2 \ln (v_k \eta_k)^2dg\\
&=\int (v_k \phi_k)^2 \left[ \ln ( (v_k \phi_k)^2 +
(v_k \eta_k)^2 ) -  \ln (v_k \phi_k)^2 \right] dg  \\
&\qquad + \int (v_k \eta_k)^2 \left[ \ln ( (v_k \phi_k)^2  +
(v_k \eta_k)^2 ) -  \ln (v_k \eta_k)^2 \right] dg \\
&\le C \int v^4_k \phi^2_k \eta^2_k dg \le C \int_{B(x_k, 4
l_k) -B(x_k,  l_k)} v^2_k \, dg. \eal
\ee Here we just used the uniform boundedness of $v_k$, proven in Lemma \ref{leyizhijie}.
This and (\ref{v2k=delta}) shows
\be
\lab{splitlog}
\int v^2_k \ln v^2_k dg =
\int (v_k \phi_k)^2 \ln (v_k \phi_k)^2 dg + \int (v_k \eta_k)^2 \ln (v_k \eta_k)^2 dg + o(1).
\ee

Recall that $v_k$ is a minimizer for $\lam_k \equiv \lam(g,
1, B(x_0, r_k) )$. By (\ref{lamk=1}),
\be
\lab{3.164}
e^{\frac{2}{n} (\lam_k -s_n)} = \frac{F(v_k)}{\exp (\frac{2}{n}
N(v_k)) }.
\ee By (\ref{splitF}) and (\ref{splitlog}), this implies
\be
\lab{fenF/N}
e^{\frac{2}{n} (\lam_k-s_n)} =\frac{F(v_k \phi_k) + F(v_k \eta_k)
}{\exp( \frac{2}{n} N(v_k \phi_k)) \, \exp( \frac{2}{n}
N(v_k \eta_k))} + o(1).
\ee Here we just used the fact that $\exp (\frac{2}{n}
N(v_k))$ is bounded away from zero. The reason is
\be
\lab{3.165}
\liminf_{k \to \infty} \exp (\frac{2}{n}
N(v_k)) = \liminf_{k \to \infty} e^{- \frac{2}{n} (\lam_k -s_n)} F(v_k)
\ge e^{- \frac{2}{n} (\lam -s_n)} F(v_\infty) >0,
\ee which is due to Step 1.

On the other hand, by definition of $\lam_k$, we have
\be
\lab{3.166}
F(v_k \phi_k) \ge e^{\frac{2}{n} (\lam_k-s_n)} \Vert v_k \phi_k
\Vert^2_{L^2} \exp \left( -\frac{2}{n} \ln \Vert v_k \phi_k
\Vert^2_{L^2} \right) \exp\left(\frac{2}{n} N(v_k \phi_k) / \Vert
v_k \phi_k \Vert^2_{L^2}\right);
\ee
\be
\lab{3.167}
F(v_k \eta_k) \ge e^{\frac{2}{n} (\lam_k-s_n)} \Vert v_k \eta_k
\Vert^2_{L^2} \exp \left( -\frac{2}{n} \ln \Vert v_k \eta_k
\Vert^2_{L^2} \right) \exp\left(\frac{2}{n} N(v_k \eta_k) / \Vert
v_k \eta_k \Vert^2_{L^2}\right).
\ee  Plugging  the last two expressions into (\ref{fenF/N}), we deduce
\be
\lab{3.168}
\frac{a^{-2/n}_k a_k A^{1/a_k}_k + b^{-2/n}_k b_k B^{1/b_k}_k  }{A_k B_k }  \le 1 + o(1),
\ee  where
\be
\lab{3.169}
a_k \equiv \Vert v_k \phi_k \Vert^2_{L^2}, \quad b_k \equiv \Vert v_k \eta_k \Vert^2_{L^2};
\ee
\be
\lab{3.170}
A_k \equiv \exp( \frac{2}{n} N(v_k \phi_k)), \quad B_k \equiv \exp( \frac{2}{n} N(v_k \eta_k)).
\ee Therefore
\be
\lab{3.171}
\min\{ a^{-2/n}_k, \,b^{-2/n}_k \} \,
\frac{ a_k A^{1/a_k}_k +  b_k B^{1/b_k}_k}{A_k B_k }  \le 1 + o(1),
\ee Notice that $a_k + b_k =1$. Therefore we have the Young's inequality: $
\frac{ a_k A^{1/a_k}_k +  b_k B^{1/b_k}_k }{A_k B_k } \ge 1$.
Letting $k \to \infty$ and using (\ref{splitv2}), we arrive at
\be
\lab{3.172}
\min\{ \delta^{-2/n}, \, (1-\delta)^{-2/n}\}  \le 1.
\ee This is a contradiction with the assumption that $ \delta =\Vert v_\infty \Vert_{L^2(M_\infty, g_\infty(0))}<1$.

\medskip

{\it Step 3.} Finally we prove that $v_\infty$ is a minimizer.

Using Fatou's Lemma, it is clear that
$F(v) \le \lim_{k \to \infty} F(v_k)$.  We claim that
\be
\lab{3.173}
N(v_\infty) \ge \lim_{k \to \infty} N(v_k),
\ee which is a reversed inequality comparing with that in Fatou's lemma.
Here goes the proof.  Let $C$ be a uniform upper bound for $\Vert v_k \Vert_\infty$.
Then $\ln (C/v_k )^2 \ge 0$. By Fatou's lemma
\be
\lab{3.174}
\int v^2_\infty \ln (C/v_\infty )^2 dg \le \lim_{k \to \infty} \int v^2_k \ln (C/v_k )^2 dg,
\ee  Since $\Vert v_\infty \Vert_{L^2} = \Vert v_k \Vert_{L^2} =1$, the above shows
\be
\lab{3.175}
N(v_\infty) = \int v^2_\infty \ln v_\infty^2 dg \ge \lim_{k \to \infty} \int v^2_k \ln v_k^2 dg =
\lim_{k \to \infty} N(v_k),
\ee  which is the claim.

Taking $k \to \infty$ in  (\ref{Fv/eN}), using the claim and  Fatou's lemma on $F(v_k)$, we deduce
\be
\lab{3.176}
e^{\frac{2}{n} (\lam -s_n)}= \lim_{k \to \infty} e^{\frac{2}{n} (\lam_k -s_n)}= \lim_{k \to \infty} \frac{F(v_k)}{\exp (\frac{2}{n}
N(v_k)) } \ge \frac{F(v_\infty)}{\exp (\frac{2}{n}
N(v_\infty)) }.
\ee Taking $\ln$ on both sides, we see that $v_\infty$ is a minimizer.
From here, it is straight forward to see that $v_\infty$ satisfies equation (\ref{maineq}).
\qed

Now we are ready to give
\medskip

\noindent {\bf  Proof of Theorem \ref{thnobreather}.}
\medskip

For simplicity, we use the notations $L(v, g) \equiv L(v, g, 1, M) $ and $\lam(g) \equiv \lam(g, 1, M)$ during the proof.

First we claim that $\lam(g)$ is invariant under scaling and  diffeomorphism. The proof is quite easy. But we
present it here to stress its independence on the behavior of the diffeomorphism at infinity.
Given any positive number $a$. It is clear that $L(v, g) = L(a^{-n/4} v, a g)$ and $\Vert v \Vert_{L^2(g)}
= \Vert a^{-n/4} v \Vert_{L^2(a g)}$. Hence $\lam(g)$ is invariant
under scaling.

Next, let $\psi$ be a diffeomorphism on $M$ and write $h = \psi^* g$. For any $v \in C^\infty_0(M)$, we have
\be
\lab{3.177}
\int_M ( 4  |\nabla v|^2 + R v^2) dg = \int_M ( 4  |\nabla_h (v \circ \psi^{-1}) |^2 + R (v \circ \psi^{-1})^2 )d h,
\ee
\be
\lab{3.178}
\int_M v^2 \ln v^2 dg = \int_M (v \circ \psi^{-1})^2 \ln (v \circ \psi^{-1})^2 dh.
\ee These imply $L(v, g) = L(v \circ \psi^{-1}, \psi^* g)$. Taking the infimum on both sides, we see that
$\lam(g)$ is also  invariant under diffeomorphism.

 Hence, we know from the assumption
$g(t_2) = c \psi^* g(t_1)$  that
\be
\lab{tk-tk}
 \lam(g(t_1)) - \lam(g(t_2))=0.
\ee

 According to Theorem \ref{thminimizer}, there exists a function $v_2 \in W^{1, 2}(M, g(t_2))$,
which is a minimizer for $\lam(g(t_2))$, i.e.
\be
\lab{l=lam}
L(v_2, g(t_2)) = L(v_2, g(t_2), 1,  M) = \lam(g(t_2)).
\ee Moreover, by Moser's iteration, it is known, as done in Lemma 2.3 in \cite{Z:2}, $v_2$ has Gaussian type decay at infinity.

Next, we solve the conjugate heat equation for $t<t_2$, with final value as $v^2_2$. This solution
is denoted by $u=u(x, t)$.
Write $v = \sqrt{u}$, then by Definition \ref{deflv}
\be
\lab{3.179}
L(v, g(t)) = - N(v) + \frac{n}{2} \ln F(v) + s_n,
\ee  where, due to $v = \sqrt{u}$,
\be
\lab{3.180}
N(v) = \int_M u \ln u \, dg(t); \quad F(v) = \int_M (\frac{|\nabla u|^2}u + R u) dg(t) =
\int_M (4 |\nabla v|^2 + R v^2) dg(t).
\ee  According to Perelman \cite{P:1} Section 1, $\frac{d}{dt} N(v) = F(v)$ and
\be
\lab{3.181}
\frac{d}{dt} F(v) = 2 \int_M | Ric - Hess (\ln u)|^2 u dg(t).
\ee We mention that although Perelman only proved the formulas for compact manifolds, but his proof also works for noncompact manifolds with bounded geometry when the functions involved have sufficiently fast decay such as the Gaussian function. See \cite{chowetc3} Chapter 19 and
 \cite{cty2007}
e.g. for a detailed computation. In our case, the function $v$ has Gaussian type decay at each time level just like the final value $v(t_2)$ does. Hence
\be
\lab{d/dtlv=}
\frac{d}{dt} L(v, g(t)) = \left( n \int_M | Ric - Hess (\ln u)|^2 u dg(t) - F^2(v) \right) \, F^{-1}(v).
\ee Following Perelman's computation,
\be
\lab{3.182}
| Ric - Hess (\ln u)|^2  \ge \left| Ric - Hess (\ln u) - \frac{1}n ( R - \Delta \ln u) g \right|^2
+ \frac{1}n ( R - \Delta \ln u)^2;
\ee Using the relation $F(v) = \int_M  (R - \Delta \ln u) u \, dg(t)$, we deduce
\be
\lab{dldt>}
\frac{d}{dt} L(\sqrt{u}, g(t)) \ge \frac{Q(u)}{F(v)}  \ge 0
\ee where
\be
\lab{qu}
\al
Q(u)(t) &= n \int_M | Ric - Hess (\ln u) - \frac{1}{n} ( R - \Delta \ln u) g |^2
u dg(t) \\
&\qquad
+ \int_M ( R - \Delta \ln u)^2 u \, dg(t)  -  \left( \int_M  (R - \Delta \ln u) u \, dg(t) \right)^2;\\
F(v)&=F(v)(t)=F(\sqrt{u})(t) = \int_M (\frac{|\nabla u|^2}u + R u) dg(t).
\eal
\ee Observe that $\sqrt{u(\cdot, t_2)} = v_2(\cdot)$ by definition. So by (\ref{l=lam}) we deduce
\be
\lab{3.183}
\al
\int^{t_2}_{t_1} \frac{d}{dt}  L(\sqrt{u}, g(t)) dt &=
L(\sqrt{u(\cdot, t_2)}, g(t_2)) - L(\sqrt{u(\cdot, t_1)}, g(t_1))\\
& \le \lam(g(t_2)) - \lam(g(t_1))=0.
\eal
\ee The last line is due to  (\ref{tk-tk}). By (\ref{dldt>}), we then have
\be
\lab{q/f}
F^{-1}(v) Q(u) =0.
\ee By (\ref{qu}), this shows that $(R - \Delta \ln u)(\cdot, t)=l(t)$, where $l=l(t)$ is a function of $t$ only.
Also
\be
\lab{3.184}
Ric - Hess (\ln u) - \frac{1}{n} l(t) g=0.
\ee Therefore, $(M, g(t))$ is a gradient Ricci soliton.
\qed

\medskip

\medskip
{\bf Acknowledgment.}  We wish to thank Professors Xiaodong Cao, Xiouxiong Chen, Zhuoran Du, Changfeng Gui and Kefeng Liu for very helpful suggestions. Thanks also go to the referees who
gave many useful input to the motivation, presentations, proofs, references in the paper.
Part of the work was done when he was a
visiting professor at Nanjing University under a Siyuan Foundation
grant, the support of which is gratefully acknowledged.

\bigskip

\noindent e-mail:  qizhang@math.ucr.edu

\enddocument